\documentclass{amsart}

\usepackage{amsfonts}
\usepackage{amsmath}
\usepackage{amsthm}
\usepackage{amssymb}
\usepackage[english]{babel}
\usepackage{graphics}
\usepackage{latexsym}
\usepackage{longtable} 
\usepackage{mathrsfs} 
\usepackage{graphicx}
\usepackage{wrapfig} 
\usepackage[pdftex,colorlinks=true,linkcolor=blue,citecolor=magenta]{hyperref}
\usepackage{enumitem}
\usepackage{esint}
\usepackage{mathcomp}
\usepackage{stmaryrd}
\usepackage{wasysym}

\newtheorem{theorem}{Theorem}
\newtheorem{corollary}[theorem]{Corollary}
\newtheorem{lemma}[theorem]{Lemma}
\newtheorem{proposition}[theorem]{Proposition}

\newtheorem{claim}[theorem]{Claim}
\newtheorem{example}[theorem]{Example}

\theoremstyle{definition}
\newtheorem{definition}[theorem]{Definition}
\newtheorem{remark}[theorem]{Remark}

\renewcommand{\d}{\mathrm{d}}

\newcommand{\mL}{\mathcal{L}}
\newcommand{\mH}{\mathcal{H}}

\newcommand{\mM}{\mathcal{M}}

\newcommand{\argmin}{\mathrm{argmin}}

\newcommand{\D}{\mathrm{D}}

\newcommand{\A}{\textbf{A}}

\newcommand{\K}{\mathrm{K}}

\newcommand{\R}{\mathbb{R}}

\newcommand{\N}{\mathbb{N}}

\newcommand{\mS}{\mathbb{S}}
\newcommand{\mB}{\mathbb{B}}

\newcommand{\noi}{\noindent}
\newcommand{\ms}{\medskip}

\newcommand{\al}{\alpha}
\newcommand{\be}{\beta}
\newcommand{\ga}{\gamma}
\newcommand{\Ga}{\Gamma}
\newcommand{\de}{\delta}
\newcommand{\De}{\Delta}
\newcommand{\e}{\varepsilon}
\newcommand{\si}{\sigma}

\newcommand{\la}{\lambda}
\newcommand{\La}{\Lambda}
\newcommand{\ka}{\kappa}
\newcommand{\Om}{\Omega}
\newcommand{\om}{\omega}

\newcommand{\av}{-\hspace{-10.5pt}\displaystyle\int}

\newcommand{\weak }{\, -\!\!\!\!\!-\!\!\!\!\rightharpoonup}
\newcommand{\weakstar }{ \overset{\, *_{\phantom{|}}}{{\smash{\, -\!\!\!\!-\!\!\!\!\rightharpoonup}}\, } }

\newcommand{\larrow}{\longrightarrow}
\newcommand{\ot}{\otimes}

\newcommand{\LL}{\text{\LARGE$\llcorner$}}
\newcommand{\p}{\partial}
\newcommand{\sub}{\subseteq}
\newcommand{\set}{\setminus}
\newcommand{\by}{\times}

\DeclareMathOperator{\sgn}{sgn} 
\newcommand{\ess}{\mathrm{ess}}
\newcommand{\diam}{\mathrm{diam}}
\newcommand{\dist}{\mathrm{dist}}

\renewcommand{\div}{\mathrm{div}}
\newcommand{\supp}{\mathrm{supp}}

\newcommand{\bt}{\begin{theorem}}\newcommand{\et}{\end{theorem}}
\newcommand{\bd}{\begin{definition}}\newcommand{\ed}{\end{definition}}
\newcommand{\bl}{\begin{lemma}}\newcommand{\el}{\end{lemma}}
\newcommand{\beq}{\begin{equation}}\newcommand{\eeq}{\end{equation}}
\newcommand{\bc}{\begin{claim}}\newcommand{\ec}{\end{claim}}
\newcommand{\bex}{\begin{example}}\newcommand{\eex}{\end{example}}
\newcommand{\bcor}{\begin{corollary}}\newcommand{\ecor}{\end{corollary}}
\newcommand{\bp}{\begin{proof}}\newcommand{\ep}{\end{proof}}

\newcommand{\BPL}{\medskip \noindent \textbf{Proof of Lemma} }

\newcommand{\BPCOR}{\medskip \noindent \textbf{Proof of Corollary} }
\newcommand{\BPP}{\medskip \noindent \textbf{Proof of Proposition} }

\numberwithin{equation}{section}

\begin{document}

\title[Vectorial $\infty$-eigenvalue nonlinear problems for $L^\infty$ functionals]{Generalised vectorial $\infty$-eigenvalue nonlinear problems for $L^\infty$ functionals}
 
\author{Nikos Katzourakis}

\address{Department of Mathematics and Statistics, University of Reading, Whiteknights Campus, Pepper Lane, Reading RG6 6AX, United Kingdom}

\email{n.katzourakis@reading.ac.uk}

  
\subjclass[2010]{35D30, 35D40, 35J47, 35J92, 35J70, 35J99, 35P30}

\date{}

\keywords{$\infty$-Eigenvalue problem; nonlinear eigenvalue problems; $\infty$-Laplacian; $L^\infty$ functionals; Absolute minimisers; Calculus of Variations in $L^\infty$; Lagrange Multipliers.}

\begin{abstract} Let $\Omega \Subset \mathbb R^n$, $f \in C^1(\mathbb R^{N\times n})$ and $g\in C^1(\mathbb R^N)$, where $N,n \in \mathbb N$. We study the minimisation problem of finding $u \in W^{1,\infty}_0(\Omega;\mathbb R^N)$ that satisfies
\[
\big\| f(\mathrm D u) \big\|_{L^\infty(\Omega)} \! = \inf \Big\{\big\| f(\mathrm D v) \big\|_{L^\infty(\Omega)} \! : \ v \! \in W^{1,\infty}_0(\Omega;\mathbb R^N), \, \| g(v) \|_{L^\infty(\Omega)}\! =1\Big\},
\]
under natural assumptions on $f,g$. This includes the $\infty$-eigenvalue problem as a special case. Herein we prove the existence of a minimiser $u_\infty$ with extra properties, derived as the limit of minimisers of approximating constrained $L^p$  problems as $p\to \infty$. A central contribution and novelty of this work is that $u_\infty$ is shown to solve a divergence PDE with measure coefficients, whose leading term is a divergence counterpart equation of the non-divergence $\infty$-Laplacian. Our results are new even in the scalar case of the $\infty$-eigenvalue problem.

\end{abstract}

\maketitle

\tableofcontents

\section{Introduction and main results}   \label{Section1}

Let $n,N \in \N$ be integers and let also $\Om \Subset \R^n$ be a bounded open set with Lipschitz boundary. In this paper we study the following variational problem: find a minimising map $u_\infty : \overline{\Om} \larrow \R^N$ in the space $W^{1,\infty}_0(\Om;\R^N)$ that solves
\beq
\label{1.1}
\big\| f(\D u_\infty) \big\|_{L^\infty(\Om)}  = \, \inf \Big\{\big\| f(\D v) \big\|_{L^\infty(\Om)}  : \ v  \in W^{1,\infty}_0(\Om;\R^N), \ \| g(v) \|_{L^\infty(\Om)} =1\Big\}.
\eeq
We are also interested in studying the structure of such constrained minimisers, and also in deriving appropriate PDEs that they satisfy as necessary conditions. The functions $f: \R^{N\by n} \larrow \R$ and $g :\R^N \larrow \R$ will be assumed to satisfy certain natural hypotheses. We note that our general notation in \eqref{1.1} and subsequently will be either self-explanatory, or a convex combination of otherwise standard symbolisations (as e.g.\ in \cite{D,FL,GM,KV}).  

The $L^p$ counterpart of \eqref{1.1}, and especially the $L^2$ case, is textbook material in the Calculus of Variations (see e.g.\ \cite{E}). Especially, the case of $f=|\cdot|$ and $g=|\cdot|$, the corresponding Euclidean norms on $\R^{N\by n}$ and $\R^N$ respectively, is known as the $p$-\emph{eigenvalue problem}, or the eigenvalue problem for the $p$-Laplacian. The $L^\infty$ case we study herein, particularly in the vectorial case of $N\geq 2$, is completely new and, despite its importance, has not been considered before. When $N=1$ with $f$ the Euclidean norm on $\R^n$ and $g$ the absolute value on $\R$, \eqref{1.1} reduces to the (scalar) $\infty$-\emph{eigenvalue problem}, whose study goes back to the seminal work of Juutinen-Lindqvist-Manfredi in \cite{JLM}. Since then, there has been a considerable interest on this problem and also on relevant ones (see for instance Juutinen-Lindqvist \cite{JL} and Bhattacharya-Marazzi \cite{BM}), as well as on related problems with constraints (see e.g.\ Aronsson-Barron \cite{AB} and Barron-Jensen \cite{BJ}). However, all these works are restricted to either the $1$-dimensional case of $n=1$, or to the scalar case of $N=1$. Crucially, these approaches rely essentially on the (scalar) Aronsson-Euler operator and on the theory of Viscosity Solutions for nonlinear PDE, both of which are not available when $N\geq 2$ (for a general introduction to the field we refer to the lecture notes \cite{C,K0}). Indeed, virtually all works on scalar-valued Calculus of Variations in $L^\infty$ rely in some way on viscosity solutions and on the comparison principle. For various interesting works, some of which are relevant to applications and some to the $\infty$-eigenvalue problem, we refer to \cite{BBJ, BJW1, BN, BBD, CDJ, CDP, GNP, HSY, KP, MWZ, PZ, RZ}. 

Vectorial and higher order variational problems involving constraints have only very recently started being explored (see \cite{K3, K4} and also \cite{CKM}). In either case, the vectorial nature of the problem requires novel methods which are not based neither on viscosity solutions, nor on the Aronsson-Euler equation. Let us also note that the mere existence of a minimiser $u_\infty$ to \eqref{1.1} is a relatively simple matter by applying the Direct Method of the Calculus of Variations and weak* lower-semicontinuity arguments for supremal functionals (under the appropriate quasi-convexity assumptions for $f$, see e.g.\ \cite{BJW2}). However, if one wishes to derive additional information on these (generally non-unique) minimisers, for example derive a necessary PDE system they satisfy, then this is a far less trivial matter. Among other complications arising, the $L^\infty$ norm is non-differentiable, non-strictly convex and lacks the property of being $\si$-additive with respect to the domain argument. 

The customary starting point for $L^\infty$ variational problems, which we also employ herein as well, is to obtain minimisers as limits of respective $L^p$ approximating variational problems as $p\to\infty$. Although this is perhaps not an intrinsic $L^\infty$ approach, it typically bears significant fruit as these \emph{special} minimisers which are obtained as limits of $L^p$ minimisers always carry a finer structure. The central novelty of this work though regards the necessary PDE conditions that such minimisers to \eqref{1.1} satisfy. Our approach in studying \eqref{1.1} is inspired by a recent development in higher order $L^\infty$ problems by the author and Moser in \cite{KM} (related earlier observations were made in a geometric higher order context in \cite{MS}). The main idea is that, if one rescales the $L^p$ Euler-Lagrange equations in a different way from the customary one used to derive the Aronsson-Euler equations as $p\to\infty$, then it is possible to derive a divergence structure PDE for the $L^\infty$ minimiser. This should be understood as a divergence form counterpart of the non-divergence Aronsson-Euler equations. However, this forcing of divergence structure to the essentially non-divergence $L^\infty$ equations has a price to be paid: the equations involve measure coefficients arising as auxiliary variables. This is somewhat reminiscent of the way that eigenvalues appear as parameters in an eigenvalue problem.

The observation that a divergence PDE with measure coefficients can be derived in $L^\infty$ has also been made earlier in a work by Evans and Yu \cite{EY}, and before that had also been conjectured by Aronsson himself, the founder of Calculus of Variations in $L^\infty$, in the unpublished note \cite{A}. However, it was not pursued further, possibly because in the case studied in \cite{EY} as a divergence form counterpart of the $\infty$-Laplacian, these measures could be highly degenerate and supported only on the boundary. This renders the PDE trivial in the interior of the domain. Nevertheless, in the higher order case employed in \cite{KM} involving the $\infty$-Bilaplacian, as well as in the case of constrained problems studied herein, these PDE systems with measure coefficients provide non-trivial information. In particular, in \cite{KM} these measure were absolutely continuous and in fact given by harmonic functions and, for \eqref{1.1}, the boundary is a nullset and these measures are supported on ``large" contact sets in the interior. 

Now we present our main results. To this end, we will utilise the following hypotheses regarding the functions $f: \R^{N\by n} \larrow \R$ and $g :\R^N \larrow \R$:
\beq
\label{1.2}
\left\{
\begin{array}{ll}
(a) & f \in C^1(\R^{N \by n}) , \ms
\\
(b) & \text{$f$ is (Morrey) quasiconvex on $\R^{N \by n}$}, \ms
\\
(c) & \text{exist $C_1,C_2>0$ with $C_1 \leq C_2$ such that} \ms
\\
& \   0  \,<  \,C_1 f(X)  \,\leq  \,\p f(X) : X  \,\leq  \,C_2 f(X), \ms
\\
& \text{for all }X \in \R^{N\by n} \set\{0\},\ms
\\
(d) & \text{exist $C_3,...,C_6>0$, $\al>1$ and $\be \leq (\al-1)/\al$ such that} \ms
\\
&  -C_3 + C_4|X|^\al  \,\leq  \,f(X) \, \leq \, C_5 |X|^\al + \,C_6,  \ms
\\
& \ \ \ \ \ \ \ \ \  \ |\p f(X)|  \,\leq  \,C_5 f(X)^\be + \,C_6,  \ms
\\
& \text{for all }X \in \R^{N\by n}, 
\end{array}
\right.
\eeq 
and
\beq
\label{1.3}
\left\{
\begin{array}{ll}
(a) & g \in C^1(\R^N), \ms
\\
(b) & \text{$g$ is coercive, i.e\ for any $\eta \in\R^N \set\{0\}$ we have} \ms
\\
& \ \ \  \  \ \ \ \ \ \ \ \ \ \ \ \ \underset{t\to\infty}{\lim} \, g(t \eta) \, = \, \infty, \ms
\\
(c) & \text{exist $C_7,C_8>0$ with $C_7\leq C_8$ such that} \ms
\\
& \ \ \ \  0 \,< \, C_7 g(\eta)  \, \leq  \,\p g(\eta) \cdot \eta  \,\leq  \,C_8 g(\eta), \ms
\\
& \text{for all }\eta \in \R^{N}\set\{0\}.
\end{array}
\right.
\eeq 
In the above, $\p f$ and $\p g$ denote the corresponding derivatives of $f$ and $g$ respectively, whilst ``$:$" and ``$\cdot$" symbolise the Euclidean inner products on $\R^{N\by n}$ and $\R^N$ respectively. By ``(Morrey) quasiconvexity" we mean the usual concept of quasiconvexity for integral functionals as e.g.\ presented in \cite{D}, not the ``$L^\infty$ quasiconvexity notions" of Barron-Jensen-Wang \cite{BJW2}. Under these hypotheses, Theorem \ref{theorem1} that follows is our first main result.

\bt \label{theorem1} Suppose that \eqref{1.2} and \eqref{1.3} hold. Then, the following are true:

\ms

\noi {\rm (A)} The problem \eqref{1.1} has a solution $u_\infty \in W^{1,\infty}_0(\Om;\R^N)$.

\ms

\noi {\rm (B)} There exist Radon measures
\[
M_\infty \in \mM(\overline{\Om};\R^{N\by n}), \ \ \ \nu_\infty \in \mM(\overline{\Om})
\]
and a number $\La_\infty \geq 0$ such that
\beq
\label{1.5}
\left\{
\begin{array}{rl}
\div ( M_\infty ) \, +\, {\La_\infty} \p g(u_\infty) \nu_\infty \,=\,0, & \text{ in }\Om,
\\
u_\infty \,=\,0, & \text{ on }\p\Om,
\end{array}
\right.
\eeq
weakly in $(C^1_0(\overline{\Om};\R^N))^*$, with
\beq
\label{1.6}
\La_\infty \, = \, \big\| f(\D u_\infty) \big\|_{L^\infty(\Om)} >\,0. 
\eeq
Further, we have the uniform bounds
\beq
\label{1.7}
\left\{ \ \ 
\begin{split}
& \La_\infty \, \geq \, \left( \! \max\left\{ \frac{C_4^{1/\al}}{ \diam(\Om) \| \p g\|_{L^\infty(\{0\leq g\leq 1\})}} -C_3^{1/\al}\,,\, 0\right\} \right)^{\!\!\al},
\\
& \La_\infty \, \leq\, C_5 \left(\! \frac{\big(g( \cdot\, \eta)\big)^{-1}(1)}{R_\Om}\right)^{\!\!\al} +\, C_6,
\end{split}
\right.
\eeq
where $R_\Om$ is the radius of the largest open ball in $\Om$, and $\big(g( \cdot\, \eta )\big)^{-1}$ is the inverse of $s \mapsto g(s \eta )$ on $[0,\infty)$ for any fixed $\eta \in \R^N$ with $|\eta|=1$ (well-defined by \eqref{1.3}).
\ms

{\rm (C)} The quadruple $\big(u_\infty, \La_\infty, M_\infty, \nu_\infty \big)$ satisfies the next approximation properties: there exists a sequence $(p_j)_1^\infty \sub (n/\al,\infty)$ with $p_j \to \infty$ as $j\to \infty$ and for any such $p$, a quadruple 
\[
\big(u_p,\La_p,M_p,\nu_p \big) \ \in \ W^{1,\al p}_0(\Om;\R^N) \by (0,\infty) \by \mM(\overline{\Om};\R^{N\by n}) \by \mM(\overline{\Om})
\]
such that
\beq
\label{1.8}
\left\{
\begin{array}{ll}
u_p \larrow u_\infty  & \text{ in }C^\ga(\overline{\Om};\R^N), \text{ for all }\ga \in(0,1),\ms
\\
\D u_p \weak \D u_\infty  & \text{ in }L^q(\Om;\R^{N\by n}), \text{ for all }q \in(1,\infty),\ms
\\
\La_p \larrow \La_\infty  & \text{ in }[0,\infty),\ms
\\
M_p \weakstar M_\infty  & \text{ in }\mM(\overline{\Om};\R^{N\by n}), \ms
\\
\nu_p \weakstar \nu_\infty  & \text{ in }\mM(\overline{\Om}), \ms
\end{array}
\right.
\eeq
as $p\to \infty$ along $(p_j)_1^\infty$. Further, $u_p$ solves the minimisation problem
\beq
\label{1.9}
\big\| f(\D u_p) \big\|_{L^p(\Om)}  = \, \inf \Big\{\big\| f(\D v) \big\|_{L^p(\Om)}  : \ v  \in W^{1,\al p}_0(\Om;\R^N), \ \| g(v) \|_{L^p(\Om)} =1\Big\}
\eeq
and $(u_p,\La_p)$ solves weakly the divergence PDE system
\beq
\label{1.10}
\left\{
\begin{array}{rl}
\div \Big( f(\D u_p)^{p-1}\p f(\D u_p)\Big)  + \, (\La_p)^p \,g(u_p)^{p-1}\p g(u_p) \,=\,0, & \text{ in }\Om,
\\
u_p \,=\,0, & \text{ on }\p\Om.
\end{array}
\right.
\eeq
Finally, the measures $M_p,\nu_p$ are given by
\beq
\label{1.11}
\left\{
\begin{split}
M_p :=& \frac{1}{\mL^n(\Om)} \left(\!\frac{f(\D u_p)}{\La_p}\right)^{\!p-1}\p f(\D u_p) \, \mL^n \LL_{\Om},
\\
\nu_p :=&  \frac{1}{\mL^n(\Om)}  g(u_p)^{p-1} \mL^n \LL_\Om .
\end{split}
\right.
\eeq
\et
We note that the meaning of the satisfaction of \eqref{1.5} weakly in $(C^1_0(\overline{\Om};\R^N))^*$ is as follows: for any test map $\phi \in C^1_0(\overline{\Om};\R^N)$, we have
\beq
\label{1.12}
\int_{\overline{\Om}} \D \phi : \d M_\infty \, = {\La_\infty} \int_{\overline{\Om}}  \p g(u_\infty) \cdot \phi
\, \d \nu_\infty .
\eeq

In general, the matrix-valued measure {\it $M_\infty$ in Theorem \ref{theorem1} implicitly depends on $\D u_\infty$}, but it is not clear how exactly, especially in the case that the derivative $\p f$ is nonlinear. Notwithstanding, if we strengthen our assumptions to include that $f$ is quadratic, together with a more stringent condition on $g$, we obtain the stronger result of Theorem \ref{theorem2} that follows. Theorem \ref{theorem2} allows to characterise \emph{$M_\infty$ as a linear function of a (special) Borel measurable representative $\D u_\infty^\star$} of the gradient of the minimiser. Additionally, we obtain much more detailed information on the structure of the measure coefficients arising. Our additional hypothesis required is the following.
\beq
\label{1.4}
\left\{
\begin{array}{l}
\text{Exists $ \A \in \R^{N\by n } \ot \R^{N\by n }$ such that, for all $X\in\R^{N\by n}$} \ms
\\
\hspace{80pt} f(X) = \A : X \ot X, \ms
\\
\text{and the assumption \eqref{1.3}(c) is satisfied with $C_7=C_8$.}
\end{array}
\right.
\eeq 
Under \eqref{1.4}, we may establish our second principal result, which is given right below.

\bt \label{theorem2} Suppose that \eqref{1.2}, \eqref{1.3} and \eqref{1.4} are satisfied. Then, in addition to the conclusions of Theorem \ref{theorem1}, the following are true:

\smallskip

\noi {\rm (A)} There exists a Radon measure $\mu_\infty \in \mM(\overline{\Om})$ and a Borel measurable mapping $\D u_\infty^\star : \overline{\Om} \larrow \R^{N\by n}$, which is a version of $\D u_\infty \in L^\infty(\Om;\R^{N\by n})$, such that
\beq
\label{1.13}
M_\infty = \, \p f(\D u_\infty^\star) \, \mu_\infty.
\eeq
Hence, $u_\infty$ solves the divergence PDE system
\beq
\label{1.14}
\left\{
\begin{array}{rl}
\div \big( \p f(\D u_\infty^\star) \, \mu_\infty\big) \, + \, {\La_\infty} \p g(u_\infty)  \nu_\infty \,=\,0, & \text{ in }\Om,
\\
u_\infty \,=\, 0, & \text{ on }\p\Om,
\end{array}
\right.
\eeq
weakly in $(C^1_0(\overline{\Om};\R^N))^*$.

\smallskip

\noi {\rm (B)} The map $\D u_\infty^\star$ can be represented as
\beq
\label{1.15}
\D u_\infty^\star (x) \, = \left\{
\begin{array}{ll}
\underset{k \to \infty}{\lim}\, \D v_{j_k} (x), & \text{if the limit exists},
\\
0, & \text{otherwise},
\end{array}
\right.
\eeq
where $(v_j)_1^\infty \sub C^1_0(\overline{\Om};\R^N)$ is any sequence satisfying that \[
\left\{ \ \ \ 
\begin{split}
\lim_{j\to\infty} \|v_j - u_\infty\|_{(W^{1,1}_0 \cap L^\infty)(\Om)} &= 0,
\\
 \limsup_{j\to\infty}\|f(\D v_j)\|_{L^\infty(\Om)} \leq \ &  \La_\infty   
\end{split}
\right.
\]
(one such is constructed in the proof). Further, $\mu_\infty$ can be approximated by
\beq
\label{1.16}
\mu_p :=  \frac{1}{\mL^n(\Om)} \left(\! \frac{f(\D u_p)}{\La_p}\right)^{\!p-1} \mL^n \LL_{\Om} ,
\eeq
as $p_j \to \infty$, in the sense that $\mu_p \weakstar \mu_\infty$ in $\mM(\overline{\Om})$. 

\smallskip

\noi {\rm (C)} The measures $\mu_\infty,\nu_\infty$ concentrate whereon $f(\D u^\star_\infty)$ and $g(u_\infty)$ are respectively maximised over $\overline{\Om}$. Additionally,
\beq
\label{1.17}
\left\{\ \ \
\begin{split}
\nu_\infty(\overline{\Om}) \, & = \, \nu_\infty \big(\big\{g(u_\infty)= 1 \big\} \big) \,=\, 1, \phantom{\Big]}
\\
\mu_\infty(\overline{\Om}) \, &  =\, \mu_\infty \left(\big\{f(\D u^\star_\infty) = \La_\infty \big\} \right) \, =\, C_8C_1^{-1} .
\end{split}
\right.
\eeq
Finally, the boundary $\p\Om$ is a nullset for both measures $\D u_\infty^\star \mu_\infty$ and $\nu_\infty$: 
\beq
\label{1.18}
\big\| \D u_\infty^\star \, \mu_\infty \big\|(\p\Om) \,= \, \nu_\infty(\p\Om) \, =\, 0.
\eeq
\et
Let us note that our results are new even in the special scalar case of the $\infty$-eigenvalue problem, as we provide a new divergence PDE describing $\infty$-eigenvectors, which is an alternative to the fully nonlinear PDE involving the $\infty$-Laplacian derived in \cite{JLM}. 

It is worth noting that our approach herein allows to bypass the need for theories of generalised solutions for the Aronsson-Euler systems arising in $L^\infty$, which are non-divergence, highly degenerate and with discontinuous coefficients. Such theories require different new ideas and some heavy machinery which depart from standard methods involving Viscosity Solutions (see e.g.\ \cite{CKP, K1} for work in this direction). Even though it is not expected that the Aronsson-Euler systems can become redundant in general, in this particular case we can indeed bypass them.

Now we describe the content and the organisation of this paper. The table of contents gives an idea regarding the order of presentation. The proofs of Theorems \ref{theorem1} and \ref{theorem2} do not appear explicitly in the main text, but instead are a consequence of numerous lemmas and propositions in Sections \ref{section2}-\ref{section6}, which gradually establish all the claims made in the statements of Theorems \ref{theorem1}-\ref{theorem2}, plus some additional auxiliary information for the $p$-problems. 

In Section \ref{section2} we discuss our hypotheses showing that, albeit restrictive, they do nonetheless allow for large classes of functions $f,g$. In particular, \eqref{1.2}(c) is compatible with the possible lack of convexity. Evidently, \eqref{1.4} is considerably more restrictive, but still allows for a large class of functions $g$.

In Section \ref{section3} we prove a large part of the assertions made in Theorem \ref{theorem1}, those which include existence for \eqref{1.1} and approximation via corresponding $L^p$ constrained minimisation problems. A new ingredient here is the necessity to construct a strongly precompact class of admissible Lipschitz maps in the respective $L^p$ minimisation classes, due to the lack of homogeneity of the functionals. Let us also note that for this part of the work, we could have relaxed our (Morrey) quasiconvexity assumption to include only ``Barron-Jensen-Wang $L^\infty$ quasiconvexity" as in \cite{BJW1,BJW2} to prove existence of minimisers, by utilising Young measures or quasi-minimisers for the approximating $L^p$ problems (for which the minima need not be attained). However, this added layer of technical complexity does not offer much insight as we need to assume considerably stronger assumptions to derive the necessary PDE systems satisfied by constrained minimisers. 

In Section \ref{section4} we establish the satisfaction of the PDE system \eqref{1.5} under the weaker assumptions \eqref{1.2}-\eqref{1.3}, by utilising tools developed in  the previous sections.

In Section \ref{section5} we introduce the appropriate mollification operators required to prove Theorem \ref{theorem2}. This regularisation scheme utilises results on the geometry of (strongly) Lipschitz domain proved in Hofmann-Mitrea-Taylor \cite{HMT}, and is closely related to the regularisation schemes used in Ern-Guermond \cite{EG}. The main idea is to use the existence of smooth vector fields which are transversal to the $\mH^{n-1}$-a.e.\ defined normal vector field on $\p\Om$ to ``shrink" the function to a compactly supported one in $\Om$, before regularising by convolution.

In Section \ref{section6} we complete the proofs of Theorems \ref{theorem1} and \ref{theorem2} by establishing the satisfaction of \eqref{1.14}, utilising the results established in earlier sections. Key ingredients here are the use of an energy identity for $L^p$ constrained minimisers and the use of Hutchinson's measure function-pairs from \cite{H}, which are a convenient way to bypass the heavy use of Young measures to identify weak* limits of sequences of products of measures with functions.

Finally, in Section \ref{section7} we consider the vectorial counterpart of the $\infty$-eigenvalue problem on the ball (i.e.\ for $f,g$ being the Euclidean norms) and in this case we are able to compute explicit measures for which the divergence PDE is satisfied.

\ms

\section{Preliminaries and discussion of hypotheses} \label{section2}

We begin by noting that in this paper we are using the notation 
\[
C^1_0(\overline{\Om};\R^N) \, :=\, C^1(\overline{\Om};\R^N) \cap C_0(\Om;\R^N)
\]
to symbolise the space of continuous mappings $\overline{\Om} \larrow \R^N$ which vanish on $\p\Om$, are $C^1$ on $\Om$ and whose derivatives extend continuously as mappings $\overline{\Om}\larrow \R^{N\by n}$. Further, everywhere in this paper, for $p \in [1,\infty)$ we are tacitly using the rescaled $L^p$-norms defined as
\[
\|h\|_{L^p(\Om)}:= \left(  \frac{1}{\mL^n(\Om)}\int_\Om|h|^p\, \d\mL^n \right)^{\!\!1/p}= \, \left(\, \av_\Om|h|^p \, \d\mL^n \right)^{\!\!1/p}, 
\]
which, by virtue of H\"older's inequality, are a family of seminorms which is monotone increasing in $p\in[1,\infty)$. Regarding the notation used in \eqref{1.4}, we note that $\A$ is a fourth order tensor, seen as a matrix over the matrix space $\R^{N\by n}$ and ``$:$" therein is the corresponding Euclidean inner product in $\R^{N\by n } \ot \R^{N\by n }$:
\[
\A : \textbf{B} \, = \sum_{1\leq \al,\be \leq N}\sum_{1\leq i,j \leq n}  \, \A_{\al i \be j} \textbf{B}_{\al i \be j}.
\]
Evidently, ``$\ot$" symbolises the tensor product of matrices, namely $X \ot Y$ is the fourth order tensor with components $X_{\al i} Y_{\be j}$, when $X,Y \in \R^{N\by n }$. We continue by noting that assumption \eqref{1.2}(c) readily implies that 
\beq
\label{2.1}
f> 0 \text{ on }\R^{N\by n}\set\{0\}, \ \ \ f(0)=0\ \text{ and $f$ is radially increasing on }\R^{N\by n},
\eeq
in the sense that $t\mapsto f(tX)$ is increasing on $(0,\infty)$, for any fixed $X\in \R^{N\by n}\set\{0\}$. Similarly, assumption \eqref{1.3}(c) implies that
\beq
\label{2.2}
g> 0 \text{ on }\R^{N}\set\{0\}, \ \ \ g(0)=0\ \text{ and $g$ is radially increasing on }\R^{N}.
\eeq
Albeit restrictive, there are nonetheless numerous (non-convex) functions which satisfy 
\eqref{1.2}-\eqref{1.3}. As an illustration, in the example below we construct a large class of functions which satisfy the inequality (c) even though they are not radially symmetric. This class of course includes all quadratic functions. We argue on $\R^{N\by n}$ only, but the construction trivially applies to $\R^{N}$ as well. Further, we remark that \eqref{1.2}(c) is compatible with \eqref{1.2}(b).

\begin{example} \label{example3} Let $\bf S$ be a compact $C^1$ hypersurface on $\R^{N\by n}$ satisfying 
\[
{\bf S} \sub \mB_{3/2}(0) \set \mB_{1/2}(0)
\]
and that any half-ray $\mathrm{span}^+[E] := \{tE: t\geq 0\}$ along a unit direction $E \in \R^{N\by n}$ intersects ${\bf S}$ at exactly one point:
\[
\mH^0(\mathrm{span}^+[E] \cap {\bf S}) = 1.
\]
We further suppose that the (outwards orientated) unit normal vector field $n_{\bf S} : {\bf S} \larrow \R^N$ satisfies $n_{\bf S}(X) :X>0$, namely the angle between the normal vector at $X$ and the direction $X$ itself is acute. Let $\al>1$ and define $f$ to be the unique $\al$-homogeneous function on $\R^{N\by n}$ satisfying
\[
{\bf S} \, =\, \{f=1\}.
\]
Explicitly, $f$ can be represented as
\[
f(X):= \inf \Big\{t>0\ : \ [0,X) \cap t^{1/\al}\bf S = \emptyset \Big\},
\]
where $[0,X)$ is the straight line segment $\{tX : 0\leq t <1\}$ and $t^{1/\al}\bf S$ is the dilation of $\bf S$ by $t^{1/\al}$. Then, by our assumptions on $\bf S$ exists $\e_0 \in (0,1)$ such that
\[
\frac{1-\e_0}{2} \, \leq\, (1-\e_0)|X| \, \leq \,  n_{\bf S}(X) :X \, \leq \,  |X|\, \leq \, \frac{3}{2},
\]
for all $X\in \bf S$. Since $\p f(X) =n_{\bf S}(X)$ and also $f(X)= 1$ for $X\in \bf S$, the above inequalities yield
\[
\frac{1-\e_0}{2}f(X) \, \leq \, \p f(X):X \, \leq \, \frac{3}{2} f(X), 
\]
for all $X\in \bf S$. Finally, since 
\[
\R^{N\by n} =\bigcup_{t\geq 0}t^{1/\al}\bf S,
\]
by noting that $f$ is $\al$-homogeneous (which implies that $\p f$ is $(\al-1)$-homogeneous), the obtained inequality in fact holds on the entire space $\R^{N\by n}$.
\end{example}

Now we derive some consequences of the satisfaction of assumptions \eqref{1.2} and \eqref{1.4}. First note that if $C_7=C_8 \equiv C$ in \eqref{1.3}, then we have $Cg(\eta)=\p g(\eta)\cdot \eta$ as equality. Even though this is restrictive, there are still several functions satisfying this condition. For instance, for any (possibly non-symmetric) matrix $A \in \R^N \ot \R^N$ and $\ga>0$, the function given by
\[
g(\eta) := (A:\eta \ot \eta)^{2\ga}
\] 
satisfies the identity with $C=4\ga$, as
\[
\p g(\eta) \cdot \eta = \Big\{ 2\ga (A:\eta \ot \eta)^{2\ga-1}A:\big((.)\ot \eta + \eta \ot (.) \big) \Big\}\cdot\eta = 4\ga \, g(\eta),
\]
for all $\eta \in \R^N$.

\begin{lemma} \label{lemma4} Suppose that $f$ satisfies  \eqref{1.2} and \eqref{1.4}. Then:
\begin{enumerate}
\item[\emph{(1)}] By setting
\[
\mS \A := \frac{\A + \A^{\top_{(1,2)\leftrightarrow (3,4)}}}{2},
\]
where $\A^{\top_{(1,2)\leftrightarrow (3,4)}}$ is the adjoint of the linear map $\A : \R^{N\by n} \larrow  \R^{N\by n}$, we have the representations
\[
f(X) =  \mS \A: X \ot X =  \big|\sqrt{\mS \A}X \big|^2,
\]
for all $X \in \R^{N\by n}$. In particular, the symmetric part $\mS \A$ of the tensor $\A$ is positive definite (i.e.\ non-negative) and hence the square root $\sqrt{\mS \A} : \R^{N\by n} \larrow \R^{N\by n}$ is well defined.

\item[\emph{(2)}] We have the identity
\[
\p f(X) = \A : \big((\cdot)\ot X + X \ot (\cdot) \big),
\]
for all $X \in \R^{N\by n}$.

\item[\emph{(3)}] $f$ is non-negative and convex on $\R^{N\by n}$ with $\{f=0\}=\{0\}$.

\item[\emph{(4)}] Let $\si$ denote the spectrum of a linear operator. Then, we have the bounds
\[
\begin{split}
\big(\!\min \si ( \mS \A )\big)|X|^2 \, \leq \, f(X) \,\leq\, |\A| |X|^2,
\end{split}
\]
for all $X \in \R^{N\by n}$. 
\end{enumerate}

\end{lemma}

Note that the adjoint (transpose) operator $\A^{\top_{(1,2)\leftrightarrow (3,4)}}$ is the operator which in index form satisfies
\[
\big({\A^{\top_{(1,2)\leftrightarrow (3,4)}}}\big)_{\al i \be j} = \, \A_{\be j \al i},
\]
for all $\al,\be \in \{1,...,N\}$ and all $i,j \in \{1,...,n\}$. This renders $\mS \A \in  \R^{N\by n} \ot \R^{N\by n}$ a real self-adjoint (symmetric) fourth order tensor (seen as a matrix over matrices). 

\BPL \ref{lemma4}. (1) The identity $f(X) =  \mS \A: X \ot X$ is immediate as the rank-one tensor $X \ot X$ is symmetric, hence $\A =\mS \A$ on the symmetric subspace of $\R^{N\by n} \ot \R^{\N\by n}$. Further, by \eqref{1.4}, \eqref{1.2} holds for $\al=2$ and by \eqref{2.1} we have
\[
\mS \A: X \ot X \geq 0 \ \ \text{ and }\ \ \mS \A: X \ot X \geq C_4|X|^2-C_3,
\] 
for all $X \in \R^{N\by n}$. By the spectral theorem it follows that the symmetric tensor $\mS \A$ is positive definite (i.e.\ non-negative) and $\si(\mS \A) \sub [0,|\A|]$. Further, the estimate above implies that $\mS \A$ is actually strictly positive and $\min \si(\mS \A) >0$: indeed, if hypothetically there existed $X_0\neq 0$ with $\mS \A: X_0 \ot X_0 =0$, then for any $t>0$
\[
0=\mS \A: (tX_0) \ot (tX_0) \geq t^2 C_4|X_0|^2-C_3,
\]
which leads to a contradiction as $t\to\infty$. Finally, since $\smash{\ms \A}$ is positive and symmetric, $ \smash{\sqrt{\mS \A}}$ exists and is also a positive symmetric operator satisfying $\mS \A =\smash{\sqrt{\mS \A}\sqrt{\mS \A}}$. Therefore,
\[
\begin{split}
f(X) &=  \mS \A: X \ot X 
\\
&= \Big(\sqrt{\mS \A}\sqrt{\mS \A}^{\top_{(1,2)\leftrightarrow (3,4)}}\Big):X\ot X
\\
&= \big(\sqrt{\mS \A}X\big) : \big(\sqrt{\mS \A}X\big)
\\
&=\big|\sqrt{\mS \A}X\big|^2,
\end{split}
\]
for any $X \in \R^{N\by n}$. 
\ms

\noi (2) Follows by a direct differentiation. 

\ms

\noi (3) Follows by part (1) by noting that $f$ is the restriction of the convex quadratic form $(X,Y) \mapsto \mS\A : X \ot Y$ on the diagonal of $\R^{N\by n} \by \R^{N\by n}$.

\ms

\noi (4) Follows directly from part (1), by recalling the variational definition of the minimum eigenvalue of a symmetric operator.     
\qed
\ms

\section{Existence in $L^p$ and compactness as $p\to \infty$}  \label{section3}

In order to solve the approximating $L^p$ constrained problems, we first need the following result which establishes the existence of strongly converging (energy comparison) maps as $p\to \infty$ in the respective constrained admissible classes.

\begin{lemma} \label{lemma5}  For any $ v\in W^{1,\infty}_0(\Om;\R^N) \set\{0\}$, there exists $(t_p)_{p\in(n/\al,\infty]} \sub (0,\infty)$ with $t_p \larrow t_\infty$ as $p\to \infty$, such that
\[
\| g(t_p v) \|_{L^p(\Om)}=\, 1
\]
for all $p\in(n/\al,\infty]$. Further, if $\| g(v) \|_{L^\infty(\Om)} =1$, then $t_\infty=1$.
\end{lemma}

\BPL \ref{lemma5}. Fix $ v\in W^{1,\infty}_0(\Om;\R^N) \set\{0\}$ and set
\[
\rho_\infty(t) := \max_{x\in \overline \Om} g ( t\,v(x)), \ \ \ t\geq 0.
\]
Then, by \eqref{2.2} and \eqref{1.3} we have that $\rho(0)=0$ and also $\rho_\infty \in C([0,\infty))$. Now we show that $\rho_\infty$ is strictly increasing. Let us begin by showing first it is non-decreasing. For any $s>0$ and $\eta\in \R^N \set\{0\}$, \eqref{1.3} implies
\[
0 < C_7 \frac{g(s \eta)}{s} \,\leq\, \p g(s \eta)\cdot \eta\, =\, \frac{\d}{\d s}\big(g(s \eta) \big),
\]
which yields that $s\mapsto g(s \eta)$ is strictly increasing on $(0,\infty)$. Then, for any $x\in \overline{\Om}$ and $t>s\geq0$ we have that $g (s\,v(x)) \leq g ( t\,v(x))$ (with strict inequality if $v(x)\neq 0$), and therefore 
\[
\rho_\infty(s) \,=\, \max_{x\in \overline{\Om}} g (s\,v(x)) \, \leq \, \max_{x\in \overline{\Om}} g (t\,v(x)) \, =\, \rho_\infty(t).
\]
Hence, $\rho_\infty$ is non-decreasing. We now show that $t\mapsto \rho_\infty(t)$ is actually strictly increasing on $(0,\infty)$. Suppose for the sake of contradiction that this is not the case. Then, by the continuity of $\rho_\infty$, there exists $t_0,\e_0>0$ such that $\rho_\infty \equiv \rho_\infty(t_0)$ on the interval $(t_0-\e_0, t_0+\e_0)$. However, by Danskin's theorem \cite{Dan}, the derivative $\rho'(t_0^+)$ from the right exists, and is given by the formula
\[
\rho_\infty'(t_0^+) \, =\, \max_{x \in \Om_{t_0}} \left\{\p g\big(t_0v(x)\big) \cdot v(x)\right\}
\]
where 
\[
\Om_{t_0} \, =\, \Big\{\bar x\in  \overline{\Om}\ : \ \rho_\infty(t_0) =g\big(t_0v(\bar x)\big)\Big\}.
\]
By utilising \eqref{1.3}, the expression for $\rho_\infty'(t_0^+)$ yields
\[
\begin{split}
\rho_\infty'(t_0^+) \, &=\, \frac{1}{t_0}\max_{x \in \Om_{t_0}} \left\{\p g\big(t_0v(x)\big) \cdot t_0 v(x)\right\} 
\\
& \geq \, \frac{C_7}{t_0}\max_{x \in \Om_{t_0}} g\big(t_0v( x)\big)
\\
& = \, \frac{C_7}{t_0}\rho_\infty(t_0 )
\\
& >\, 0.
\end{split}
\]
This is a contradiction to $\rho_\infty$ being constant on $(t_0, t_0+\e_0)$, hence establishing that $\rho_\infty$ is indeed strictly increasing. Recall now that by \eqref{1.3}(b) we also have $g(s \eta) \to \infty$ as $s\to\infty$.  Thus, for any fixed $\bar x \in \Om$ for which $v(\bar x)\neq 0$, we have
\[
\lim_{t\to \infty} \rho_\infty(t) \, \geq \, \lim_{t\to \infty} g\big( t\, v(\bar x)\big) \, =\, \infty.
\]
The mean value theorem then implies that there exists a number $t_\infty >0$ such that $\rho_\infty(t_\infty)=1$, namely
\[
\| g( t_\infty v ) \|_{L^\infty(\Om)} =\, 1.
\]
Additionally, if $\| g( v ) \|_{L^\infty(\Om)}=1$ to begin with, then $t_\infty=1$. Fix now $p \in (n/\al,\infty)$ and set
\[
\rho_p(t)\, := \, \av_\Om g\big( t\, v(x) \big)^p \, \d x , \ \ \ t\geq 0.
\]
Then, by \eqref{2.2} we have $\rho_p(0)=0$ and also the strict monotonicity of $s\mapsto g(se)$ on $(0,\infty)$. Then, the monotone convergence theorem yields  that $\rho_p(t) \to \infty$ as $t\to\infty$. Therefore, there exists $t_p >0$ such that $\rho_p(t_p)=1$, namely
\[
\| g( t_p v ) \|_{L^p(\Om)} =\, 1.
\]
We now claim that $t_p \larrow t_\infty$ as $p\to \infty$. Assuming for the sake of contradiction that $t_p \not\!\!\larrow t_\infty$ as $p\to\infty$, this means there exists a sequence $(t_{p_j})_1^\infty \sub (n/\al,\infty)$ and a $t_0 \in [0,t_\infty)\cup(t_\infty,\infty]$ such that $t_{p_j} \larrow t_0$ as $j\to \infty$. Further, the sequence $(t_{p_j})_1^\infty$ can be assumed to be monotone (either increasing or decreasing). We first show that $t_0<\infty$. Indeed, if we had $t_0 =\infty$, then the sequence $(t_{p_j})_1^\infty$ must be monotone increasing. Hence, by \eqref{1.3}(b) we have $g( t_{p_j} v(x)) \nearrow \infty$ as $j\to \infty$, for any $x\in \Om \cap \{v\neq 0\}$. The monotone convergence theorem then yields the contradiction
\[
1 \, =\, \av_\Om g ( t_{p_j} v )^{p_j} \, \d \mL^n \, \nearrow \,\infty.
\]
Therefore, we must have $t_0<\infty$ and hence $t_0 \in [0,t_\infty)\cup(t_\infty,\infty)$. If now $t_{p_j} \larrow t_0$, we have $t_{p_j}v \larrow t_0v$ in $C(\overline{\Om};\R^N)$, as $j\to \infty$. By the continuity of $g$, this implies
$ g ( t_{p_j} v ) \larrow g ( t_0 v )$ in $C(\overline{\Om})$, as $j\to \infty$. Since $\|\cdot\|_{L^p(\Om)}\larrow \|\cdot\|_{L^\infty(\Om)}$ as $p\to \infty$ in the pointwise sense on $L^\infty(\Om)$, we infer that
\[ 
\begin{split}
1 \, &= \, \| g ( t_{p_j} v ) \|_{L^{p_j}(\Om)} 
\\
& = \, \| g ( t_0 v ) \|_{L^{p_j}(\Om)} \, + \, O\big(\big\| g ( t_{p_j}v ) -g ( t_0 v  ) \big\|_{L^\infty(\Om)} \big)
\\
& = \, \| g ( t_0 v ) \|_{L^{p_j}(\Om)} +  O\big(| t_{p_j} - t_0| \big)
\\
&\!\!\! \overset{j\to\infty}{\larrow} \ \| g ( t_0 v )  \|_{L^\infty(\Om)}
\\
& = \,\rho_\infty(t_0),
\end{split}
\]
which is a contradiction if $t_\infty \neq t_0$ due to the strict monotonicity of $\rho_\infty$.        \qed

\ms

Now we proceed to show existence of solution to the approximate constrained minimisation problem.

\begin{lemma} \label{lemma6} For any $p>n/\al$, the minimisation problem \eqref{1.9} has a solution $u_p \in W^{1,\al p}_0(\Om;\R^N)$.
\end{lemma}

\BPL \ref{lemma6}. Let $p\in(n/\al,\infty)$ be fixed. We begin by noting that, by virtue of Lemma \ref{lemma5}, the admissible class is non-empty: for any fixed $v_0 \in W^{1,\infty}_0(\Om;\R^N)$ with $v_0 \not\equiv 0$,  there exists $t_p>0$ such that $\| g(t_p v_0) \|_{L^p(\Om)}=1$, yielding that $t_p v_0$ is in the admissible class of \eqref{1.9}. Next, since by \eqref{1.2}(b) $f$  is (Morrey) quasiconvex in $\R^{N\by n}$, so is $f^p$, as a result of Jensen's inequality: indeed, for any fixed $X\in \R^{N\by n}$ and $\phi \in W^{1,\infty}_0(\Om;\R^N)$, we have
\[
f^p(X) \,\leq \Bigg(\, {\av}_{\!\!\!\Om} f(X+\D \phi)\,\d\mL^n \!\Bigg)^{\!\!p} \,\leq \, {\av}_{\!\!\!\Om} f^p(X+\D \phi)\,\d\mL^n .
\]
Further, by \eqref{1.2}(d) $f^p$ satisfies the estimate
\[
 -C_3(p) + C_4(p)|X|^{\al p} \, \leq \,f^p(X) \,\leq \, C_5(p) |X|^{\al p} +C_6(p),
\]
for some $p$-dependent constants $C_3(p),...,C_6(p)>0$. By standard results on quasiconvex integral functionals (see e.g.\ \cite{D}), it follows that $\|f(\D (\cdot))\|_{L^p(\Om)}$ is weakly lower semi-continuous and coercive in the space $W^{1,\al p}_0(\Om;\R^N)$. Further, by the Morrey embedding $W^{1,\al p}_0(\Om;\R^N) \sub C^{1-\frac{n}{\al p}}(\overline{\Om};\R^N)$, it follows that the admissible class is weakly closed as the functional $\|g(\cdot)\|_{L^p(\Om)}$ is weakly continuous on $W^{1,\al p}_0(\Om;\R^N)$. Hence, there exists a minimiser $u_p$ which solves \eqref{1.9}, as claimed.   \qed

\ms

Now we consider the Euler-Lagrange equations that the approximate minimiser $u_p$ satisfies. Unsurprisingly, they involve a Lagrange multiplier, which arises by the integral constraint $\|g(\cdot)\|_{L^p(\Om)}=1$.

\begin{lemma} \label{lemma7} For any $p>n/\al$, let $u_p$ be the minimiser of \eqref{1.9} given by Lemma \ref{lemma6}. Then, there exists $\la_p \in \R$ such that the pair $(u_p,\la_p) \in W^{1,\al p}_0(\Om;\R^N) \by \R $  satisfies weakly the PDE system
\[
\left\{
\begin{array}{rl}
\div \Big( f(\D u_p)^{p-1}\p f(\D u_p)\Big)  + \, \la_p \,g(u_p)^{p-1}\p g(u_p) \,=\,0, & \text{ in }\Om,
\\
u_p \,=\,0, & \text{ on }\p\Om.
\end{array}
\right.
\]
\end{lemma}

\BPL \ref{lemma7}. Follows by standard results on Lagrange multipliers in Banach spaces (see e.g.\ Zeidler \cite[Th.\ 43.D, p.\ 290]{Z}). \qed
\ms

Now we obtain additional information on the family of eigenvalues $(\la_p)_{p>n/\al}$.

\begin{lemma} \label{lemma8} For any $p>n/\al$, we set
\[
L_p \,:=\, \big\|f(\D u_p)\big\|_{L^p(\Om)}.
\]
Then, there exists $\La_p>0$ such that
\[
\la_p \,= \,(\La_p)^p
\]
and also
\[
0\,< \, \bigg(\frac{C_1}{C_8}\bigg)^{\!\!1/p} L_p \, \leq\, \La_p \,\leq\, \bigg(\frac{C_2}{C_7}\bigg)^{\!\!1/p} L_p.
\]
\end{lemma}

\BPL \ref{lemma8}. We begin by noting that $L_p>0$, namely the infimum over the admissible class in \eqref{1.9} is strictly positive, as a consequence  of the constraint and our assumptions \eqref{1.2}-\eqref{1.3} (via \eqref{2.1}-\eqref{2.2}): the only map $u \in W^{1,\al p}_0(\Om;\R^N)$ for which $\|f(\D u)\|_{L^p(\Om)}=0$ is $u_0\equiv 0$, but $u_0$ is not in the admissible class because $\|g(u_0)\|_{L^p(\Om)}=0$. Next, by testing against $u_p$ in the weak formulation of the Euler-Lagrange equations (Lemma \ref{lemma7}), we have
\[
\int_\Om f(\D u_p)^{p-1}\p f(\D u_p) :\D u_p \,\d\mL^n \,=\, \la_p \int_\Om g(u_p)^{p-1}\p g(u_p) \cdot u_p  \,\d\mL^n. 
\]
By \eqref{1.2}(c) and \eqref{1.3}(c) we have
\[
\left\{ \ \ 
\begin{split}
 C_1 f(\D u_p) \,\leq \,\p f(\D u_p) : \D u_p \,\leq\, C_2 f(\D u_p),
 \\
 C_7 g( u_p) \,\leq\, \p g( u_p) \cdot u_p \,\leq \,C_8 g(u_p), \ \ \ \
\end{split}
\right.
\]
$\mL^n$-a.e.\ on $\Om$. Hence, since $f,g\geq0$, integration gives
\[
\left\{ \ \, 
\begin{split}
 C_1 \ {\av}_{\!\!\!\Om}f^p(\D u_p) \,\d\mL^n \, \leq \ {\av}_{\!\!\!\Om} f^{p-1}(\D u_p)\p f(\D u_p) : \D u_p \,\d\mL^n \, \leq \ C_2 \ {\av}_{\!\!\!\Om} f^p(\D u_p)\,\d\mL^n,
 \\
 C_7 \ {\av}_{\!\!\!\Om}g^p( u_p) \ \d\mL^n\leq \ {\av}_{\!\!\!\Om}g^{p-1}( u_p) \p g( u_p) \cdot u_p \,\d\mL^n \, \leq \ C_8 \ {\av}_{\!\!\!\Om} g^p(u_p)\,\d\mL^n. \ \ \ \
\end{split}
\right.
\]
By recalling that
\[
{\av}_{\!\!\!\Om}f^p(\D u_p) \,\d\mL^n \,=\, (L_p)^p \, > \, 0 \ \ \ \text{ and } \ \ \ {\av}_{\!\!\!\Om} g^p(u_p)\,\d\mL^n \,=\, 1,
\]
it follows that $\la_p>0$. By virtue of the above, by defining $\La_p := (\la_p)^{1/p} >0$, we have the estimates 
\[
C_1 (L_p)^p \, \leq \, \la_pC_8 ,\ \ \ \ C_2  (L_p)^p \, \geq \, \la_p C_7,
\]
which lead to the desired inequality.
\qed

\ms

We now show the existence of solution to the problem \eqref{1.1} and the compactness of the class of $p$-pairs of eigenvectors-eigenvalues $(u_p,\La_p)_{p>n/\al}$.

\begin{proposition} \label{proposition9} There exists $(u_\infty,\La_\infty) \in W^{1,\infty}_0(\Om;\R^N) \by (0,\infty)$ such that, along a sequence $(p_j)_1^\infty$ we have
\[
\left\{
\begin{array}{ll}
u_p \larrow u_\infty  & \text{ in }C^\ga(\overline{\Om};\R^N), \text{ for all }\ga \in(0,1),\ms
\\
\D u_p \weak \D u_\infty  & \text{ in }L^q(\Om;\R^{N\by n}), \text{ for all }q \in(1,\infty),\ms
\\
\La_p \larrow \La_\infty  & \text{ in }(0,\infty).
\end{array}
\right.
\]
Further, $u_\infty$ solves the minimisation problem \eqref{1.1} and also $\La_\infty$ is given by \eqref{1.6}. Finally, $\La_\infty$ satisfies the uniform bounds \eqref{1.7}.
\end{proposition}

\BPP \ref{proposition9}. Fix $p>n/\al$, $q\leq p$ and a map $ v_0\in W^{1,\infty}_0(\Om;\R^N)$ with $v_0 \not \equiv 0$. Then, by Lemma \ref{lemma5} there exists $(t_p)_{p\in(n/\al,\infty]} \sub (0,\infty)$ with $t_p \larrow t_\infty$ as $p\to \infty$, such that $\| g(t_p v_0) \|_{L^p(\Om)}=1$ for all $p\in(n/\al,\infty]$. By H\"older's inequality and minimality, this allows to estimate
\[
\|f(\D u_p)\|_{L^q(\Om)} \, \leq\, \|f(\D u_p)\|_{L^p(\Om)} \, \leq\, \|f(\D (t_p v_0))\|_{L^p(\Om)} \, \leq\, \|f(t_p \D v_0)\|_{L^\infty(\Om)}.
\]
Hence,
\[
\sup_{q\geq p} \|f(\D u_p)\|_{L^q(\Om)} \, \leq\, \sup_{q\geq p}\|f(t_p \D v_0)\|_{L^\infty(\Om)} <\infty,
\]
and the finiteness of the last term is a consequence of the compactness of $(t_p)_{p>n/\al}$. Further, by \eqref{1.2} we have the lower bound $f^q(X) \geq C_4(q)|X|^{\al q} -C_3(q)$ for some $q$-dependent constants $C_3(q),C_4(q)>0$ and all $X\in \R^{N\by n}$. When combined with the previous estimate, this lower bound implies
\[
\sup_{q\geq p} \|\D u_p\|_{L^{\al q}(\Om)}  \, \leq  \, C(q) <\infty,
\]
for some $q$-dependent $C(q)>0$. Further, application of Poincar\'e's inequality improves the above estimate to
\[
\sup_{q\geq p} \| u_p\|_{W^{1,\al q}(\Om)}  \, \leq C(q) \,  <\infty,
\]
for a new constant $C(q)>0$. Hence, standard compactness and diagonal arguments in Sobolev spaces imply the existence of a map
\[
u_\infty \, \in \bigcap_{n/\al < p<\infty} W_0^{1,\al p}(\Om;\R^N)
\]
and (for any sequence) a subsequence such that the claimed modes of convergence hold true as $p\to \infty$ along this subsequence. Fix now any $v\in W_0^{1,\infty}(\Om;\R^N)$ satisfying $\|g(v)\|_{L^\infty(\Om)}=1$. Note that, necessarily, $v\not\equiv 0$ by virtue of \eqref{2.2}. By Lemma \ref{lemma5}, there exists $(t_p)_{p\in(n/\al,\infty)} \sub (0,\infty)$ with $t_p \larrow 1$ as $p\to \infty$, such that $\| g(t_p v) \|_{L^p(\Om)}=1$ for all $p>n/\al$. By Lemma \ref{lemma6}, the definition of $L_p$ in Lemma \ref{lemma8}, H\"older inequality and minimality, we estimate
\[
\|f(\D u_p)\|_{L^q(\Om)}  \, \leq\, L_p\, \leq\, \|f(t_p \D v)\|_{L^p(\Om)}
\]
for any such $v$. By the weak lower semi-continuity of the functional $\|f(\D (\cdot))\|_{L^q(\Om)} $ on $W_0^{1,\al q}(\Om)$, by letting $p\to \infty$ along a subsequence, we obtain for the given fixed $v\in W_0^{1,\infty}(\Om;\R^N)$ that
\[
\|f(\D u_\infty)\|_{L^q(\Om)}  \, \leq\, \liminf_{p_j\to \infty} L_p\,  \, \leq\, \limsup_{p_j\to \infty} L_p\, \leq\, \|f(\D v)\|_{L^\infty(\Om)}.
\]
In conclusion, by letting also $q\to\infty$, we deduce the energy inequality
\[
\left\{ \ \ 
\begin{split}
& \|f(\D u_\infty)\|_{L^\infty(\Om)}  \, \leq\, \liminf_{p_j\to \infty} L_p\,  \, \leq\, \limsup_{p_j\to \infty} L_p\, \leq\, \|f(\D v)\|_{L^\infty(\Om)}, 
\\
& \ \ \ \ \ \ \ \text{for any $v\in W_0^{1,\infty}(\Om;\R^N)$ with $\|g(v)\|_{L^\infty(\Om)}=1$.}
\end{split}
\right.
\]
Next, we note that $\D u_\infty \in L^\infty(\Om;\R^{N\by n})$ (and is not merely in the intersection of $L^{\al p}(\Om;\R^{N\by n})$ for $p\in (n/\al,\infty)$). This follows by \eqref{1.2}, which yields
\[
\|f(\D u_\infty)\|_{L^\infty(\Om)} \, \geq \, C_4\|\D u_\infty\|^\al_{L^\infty(\Om)} - \, C_3.
\]
As a result, by Morrey's estimate we have that $u_\infty \in W^{1,\infty}_0(\Om;\R^{N})$. Further, $u_\infty$ is in the admissible class of \eqref{1.1}, since by the uniform convergence $u_p \larrow u_\infty$ on $\overline{\Om}$ as $p\to \infty$ along a subsequence, the continuity of $g$, we have 
\[
\begin{split}
1 \, & =\, \|g(u_p) \|_{L^p(\Om)} 
\\
& =\, \|g(u_\infty) \|_{L^p(\Om)} \, + \, \big( \|g(u_p) \|_{L^p(\Om)} -\,  \|g(u_\infty) \|_{L^p(\Om)} \big)
\\
& =\, \|g(u_\infty) \|_{L^p(\Om)} \, + \, O\big( \|g(u_p)  - g(u_\infty) \|_{L^p(\Om)} \big)
\\
& =\, \|g(u_\infty) \|_{L^p(\Om)} \, + \, O\big( \|g(u_p)  - g(u_\infty) \|_{L^\infty(\Om)} \big)
\\
& \!\! \larrow \, \|g(u_\infty) \|_{L^\infty(\Om)},
\end{split}
\]
as $p\to \infty$ subsequentially. In conclusion, by the arbitrariness of the map $v$ in our energy inequality, it follows that $u_\infty$ indeed solves the minimisation problem \eqref{1.1}. 

Let now $\La_\infty$ be {\it defined} by
\[
\La_\infty\, :=\, \|f(\D u_\infty)\|_{L^\infty(\Om)}. 
\]
Note that the above definition of $\La_\infty$ as the infimum in \eqref{1.1} readily implies that it must be strictly positive. Indeed, by \eqref{2.1}-\eqref{2.2}, for any map $u_0$ in the admissible class of \eqref{1.1} for which we have $\|f(\D u_0)\|_{L^\infty(\Om)}=0$, it follows that $u_0 \equiv 0$, which is a contradiction since we must also have $\|g(u_0)\|_{L^\infty(\Om)}=0 <1$. Further, the choice $v:= u_\infty$ in our energy inequality implies that $L_p \larrow \La_\infty$ as $p\to\infty$ along a subsequence. Additionally, by Lemma \ref{lemma8} we have that $|L_p -\La_p|\larrow 0$ as $p\to \infty$, hence $\La_p \larrow \La_\infty$ as $p\to \infty$. 

To complete the proof, it remains to establish the claimed uniform bounds on $\La_\infty$. By Poincar\'e's inequality, since $g(0)=0$, we have
\[
\begin{split}
1 \,= \, \|g(u_\infty)\|_{L^\infty(\Om)} \, & \leq \, \diam(\Om) \|\D(g(u_\infty))\|_{L^\infty(\Om)}
\\
& \leq  \,\diam(\Om) \|\p g\|_{L^\infty(u_\infty(\overline{\Om}))} \|\D u_\infty\|_{L^\infty(\Om)}.
\end{split}
\]
Since $g\geq 0$ and $\|g(u_\infty)\|_{L^\infty(\Om)}=1$, this yields that $0\leq g(u_\infty) \leq 1$ everywhere on $\overline{\Om}$. Hence, we have $u(\overline{\Om}) \sub \{0\leq g\leq 1\}$, which in turn implies
\[
\begin{split}
1  \, \leq  \,\diam(\Om) \|\p g\|_{L^\infty(\{0\leq g\leq 1\})} \|\D u_\infty\|_{L^\infty(\Om)}.
\end{split}
\]
Since by \eqref{1.2} we have $|X|\leq C_4^{-1/\al}\big(f(X)+C_3\big)^{1/\al}$ for any $X\in \R^{N\by n}$, we obtain
\[
\begin{split}
1 \,\leq \, \diam(\Om) \|\p g\|_{L^\infty(\{0\leq g\leq 1\})}  \frac{ \|f(\D u_\infty)\|^{1/\al}_{L^\infty(\Om)} + \, C_3^{1/\al} }{C_4^{1/\al}}   .
\end{split}
\]
The above inequality readily leads to the claimed lower bound for $\La_\infty$ in \eqref{1.7}.

For the upper bound, consider the Lipschitz function $\dist(\cdot,\p\Om)$ of distance from the boundary $\p\Om$, extended by zero on $\R^n \set\Om$. Since $\dist(\cdot,\p\Om) \in W^{1,\infty}_0(\Om)$, for any fixed $\eta \in \R^N$ with $|\eta|=1$ we may invoke Lemma \ref{lemma5} to find $t_\infty >0$ such that $\xi_\infty :=t_\infty \dist(\cdot,\p\Om) \eta $ is in the admissible class for \eqref{1.1}, namely $\|g(\xi_\infty)\|_{L^\infty(\Om)}=1$. Since $|\D(\dist(\cdot,\p\Om))|\leq 1$ $\mL^n$-a.e.\ on $\Om$, by minimality and \eqref{1.2}, we estimate
\[
\begin{split}
\La_\infty \, & =\, \|f(\D u_\infty) \|_{L^\infty(\Om)}
\\
& \leq \, \|f(\D \xi_\infty) \|_{L^\infty(\Om)}
\\
& \leq \, \big\|f\big(t_\infty \eta  \ot \D(\dist(\cdot,\p\Om)) \big) \big\|_{L^\infty(\Om)}
\\
& \leq \, C_5 \Big(t_\infty \big\| \eta  \ot \D(\dist(\cdot,\p\Om)) \big\|_{L^\infty(\Om)}\Big)^{\!\al} +\, C_6
\\
& \leq \, C_5 (t_\infty)^\al +\, C_6.
\end{split}
\] 
Now we need an estimate for $t_\infty$. To this end, let
\[
R_\Om \,:=\, \sup \{r>0 \ | \ \exists \ x \in \Om \, : \ \mB_r(x) \sub \Om \}.
\]
Then, there exists an $\bar x \in \Om$ such that $\mB_{R_\Om}(\bar x) \sub \Om$, and this is the largest such ball. Hence, $\dist(\bar x, \p\Om) =R_\Om$. Since by construction, $\xi_\infty$ satisfies $\|g(\xi_\infty)\|_{L^\infty(\Om)}=1$, it follows that
\[
1 \, =\, \sup_\Om \, g\big(t_\infty \dist (\cdot,\p\Om) \eta  \big) \, \geq\, g\big(t_\infty \dist (\bar x,\p\Om)\eta  \big) \, =\, g(t_\infty R_\Om \eta ).
\]
By \eqref{2.2}, we have that the function $[0,\infty) \ni s \mapsto g(s\eta ) \in [0,\infty)$ is strictly increasing and onto because $g(0)=0$ and $g(s \eta )\to \infty$ as $s\to\infty$. Hence, its inverse function $(g(\cdot \, \eta ))^{-1} : [0,\infty) \larrow [0,\infty)$ is well defined. Therefore, since $g(t_\infty R_\Om \eta ) \leq 1$,  we have
\[
t_\infty R_\Om  \, =\, (g(\cdot \, \eta ))^{-1}\big( g(R_\Om \eta )\big) \, \leq \, (g(\cdot \, \eta ))^{-1}(1).
\]
Combining the previous estimates, we infer that
\[
\La_\infty \, \leq\, C_5 \bigg(\! \frac{1}{R_\Om} \big(g(\, \cdot\, \eta )\big)^{-1}(1)\bigg)^{\!\!\al} +\, C_6.
\]
The above estimate completes the claimed uniform upper bound on $\La_\infty$ in \eqref{1.7}, and therefore Proposition \ref{proposition9} ensues. \qed

\ms

\section{The divergence PDE system in $L^\infty$, part I} \label{section4}

Using the tools already developed, in this section we establish the satisfaction of the divergence PDE system \eqref{1.5} by the pair $(u_\infty,\La_\infty)$, under only assumptions \eqref{1.2}-\eqref{1.3}. Later in Section \ref{section6} we will establish satisfaction of \eqref{1.14} under the additional stronger assumption \eqref{1.4}. We begin with the limiting measures.

\begin{lemma} \label{lemma10} For any $p>n/\al+2$, consider the non-negative measures $\mu_p,\nu_p \in \mM(\overline{\Om})$ and also the matrix-valued measure $M_p \in \mM(\overline{\Om};\R^{N\by n})$, given by \eqref{1.11} and \eqref{1.16}. Then, for all large enough $p$ we have the bounds
\[
\left\{
\begin{split}
 \mu_p (\overline{\Om}) & \leq \,  \bigg(\frac{C_8}{C_1}\bigg)^{\!1-\frac{1}{p}},
\\
\| M_p \|(\overline{\Om}) & \leq \,  \bigg(\frac{C_8}{C_1}\bigg)^{\!1-\frac{1}{p}}\Big(C_5 \big(\La_\infty +1 \big)^{\be}+ C_6\Big) ,
\end{split}
\right.
\]
and, by setting $\om(p):=\big\| g(u_p)-g(u_\infty )\big\|_{L^\infty(\Om)}$, we have 
\[
\left\{ \ \ 
\begin{split}
 \frac{ 1 }{ 1+ \om(p)} \, \leq \, \nu_p (\overline{\Om}) \,  \leq \, 1 ,& 
\\
\nu_p \Big(\Big\{g(u_\infty)< \|g(u_\infty)\|_{L^\infty(\Om)} - 2\om(q) \Big\} \Big) \, &\leq\, \big(1-\om(q)\big)^{p-1},
\end{split}
\right.
\]
for $p,q \geq n/\al+2$ large enough. Consequently, there exists a further subsequence and limiting measures $\mu_\infty,\nu_\infty \in \mM(\overline{\Om})$ and $M_\infty \in \mM(\overline{\Om};\R^{N\by n})$ such that $\mu_p \weakstar \mu_\infty$ and $\nu_p \weakstar \nu_\infty$ in $\mM(\overline{\Om})$, and also $M_p \weakstar M_\infty$ in $\mM(\overline{\Om};\R^{N\by n})$, as $p\to\infty$, along this subsequence. Additionally, $\nu_\infty$ satisfies
\[
\nu_\infty(\overline{\Om}) = \nu_\infty \Big(\Big\{g(u_\infty)= \|g(u_\infty)\|_{L^\infty(\Om)} \Big\} \Big) = 1, \ \ \ \nu_\infty(\p\Om)=0.
\]
\end{lemma}

\BPL \ref{lemma10}. By \eqref{1.16} and Lemma \ref{lemma8} we estimate 
\[
\begin{split}
\mu_p (\Om)\, & = \,   {\av}_{\!\!\!\Om} \bigg(\frac{f(\D u_p)}{\La_p}\bigg)^{\!p-1} \d \mL^n
\\
& \leq  \, \frac{1}{\La_p^{p-1}} \left(\, \av_\Om  f(\D u_p)^p \, \d \mL^n \!\right)^{\!\!\frac{p-1}{p}}
\\
& =  \,   \bigg(\frac{L_p}{\La_p}\bigg)^{\!p-1}
\\
& \leq  \,  \bigg(\frac{C_8}{C_1}\bigg)^{\!1-\frac{1}{p}}.
\end{split}
\]
Further, for $p>n/\al+2$, by \eqref{1.11} and assumption \eqref{1.2}(d), we estimate
\[
\begin{split}
\| M_p \|(\overline{\Om})  \, &= \,  {\av}_{\!\!\!\Om} \bigg(\frac{f(\D u_p)}{\La_p}\bigg)^{\!p-1} \big|\p f(\D u_p) \big| \, \d \mL^n
\\
&\leq \, \frac{1}{\La_p^{p-1}}  \av_\Om  f(\D u_p)^{p-1} \Big(C_5 f(\D u_p)^\be +C_6\Big)\, \d \mL^n 
\\
&= \, \frac{1}{\La_p^{p-1}} \av_\Om  \Big(C_5 f(\D u_p)^{p-1+\be} + C_6 f(\D u_p)^{p-1}\Big)\, \d \mL^n  .
\end{split}
\]
In view of Lemma \ref{lemma8} and Proposition \ref{proposition9}, for large $p$ the previous estimate yields
\[
\begin{split}
\| M_p \|(\overline{\Om}) \, &\leq \, \frac{C_5}{\La_p^{p-1}} \! \left( \, \av_\Om f(\D u_p)^p \, \d \mL^n \right)^{\! \frac{p-1+\be}{p}}  +\, \frac{C_6}{\La_p^{p-1}} \! \left( \, \av_\Om f(\D u_p)^p \, \d \mL^n \right)^{\! \frac{p-1}{p}}    
\\
&= \, C_5 \frac{L_p^{p-1+\be}}{\La_p^{p-1}} \,  +\, C_6 \frac{L_p^{p-1}}{\La_p^{p-1}} 
\\
&= \, \bigg(\frac{L_p}{\La_p}\bigg)^{p-1}\big(C_5 L_p^\be + C_6 \big)
\\
& \leq\, \bigg(\frac{C_8}{C_1}\bigg)^{\!1-\frac{1}{p}}\Big(C_5 (\La_\infty+1)^\be + C_6 \Big),
\end{split}
\]
as claimed. Similarly, by \eqref{1.11} and the fact that $\|g(u_p)\|_{L^p(\Om)}=1$, we may estimate
\[
\begin{split}
 \nu_p (\overline{\Om}) \,& =  \,   {\av}_{\!\!\!\Om}  g(u_p)^{p-1} \,\d \mL^n
\, \leq    \left(\, \av_\Om  g(u_p)^p \, \d \mL^n \right)^{\!\!\frac{p-1}{p}}
\, =  \, 1 .
\end{split}
\]
Next, by setting $\om(p):=\| g(u_p)-g(u_\infty )\|_{L^\infty(\Om)} $ and noting that $\om(p)\to0$ as $p\to \infty$ (along a subsequence) due to the uniform convergence $g(u_p)\larrow g(u_\infty )$ on $\overline{\Om}$, we estimate
\[
\begin{split}
1 \, & = \, \av_\Om  g(u_p)^{p} \,\d \mL^n
\\
& = \,  \av_\Om  g(u_p)^{p-1} g(u_p) \, \d \mL^n  
\\
& \leq \,  \av_\Om  g(u_p)^{p-1} \Big(g(u_\infty) + \| g(u_p)-g(u_\infty )\|_{L^\infty(\Om)}  \Big) \, \d \mL^n  
\\
& \leq \,  \av_\Om  g(u_p)^{p-1}  \big( 1 + \om(p)\big) \, \d \mL^n , 
\\
\end{split}
\]
where in the last step we also used that $g(u_\infty)\leq 1$ on $\Om$, as a consequence of the continuity and the non-negativity of $g(u_\infty)$ and the constraint $\|g(u_\infty)\|_{L^\infty(\Om)}=1$. Hence, the above inequality shows that
\[
1 \, \leq\, \nu_p(\overline{\Om})\big( 1 + \om(p)\big).
\]
Finally, let us fix $q\geq n/\al+2$ and define the open set
\[
A_q \,:=\, \Big\{g(u_\infty)< \|g(u_\infty)\|_{L^\infty(\Om)} - 2\om(q) \Big\}.
\]
We note that the desired remaining estimate for $\nu_p$ can be deduced by the more general result \cite[Proposition 7]{K2}, but for the sake of completeness we provide a self-contained simpler proof. For $p$ and $\e$ as above, we estimate
\[
\begin{split}
\nu_p(A_q) & = \,  \frac{1}{\mL^n(\Om)} \int_{A_q} g(u_p)^{p-1} \,\d \mL^n
\\
& \leq \,  \frac{1}{\mL^n(\Om)} \int_{A_q} \Big(g(u_\infty) + \| g(u_p)-g(u_\infty )\|_{L^\infty(\Om)}  \Big)^{p-1} \,\d \mL^n
\\
& \leq \, \frac{1}{\mL^n(\Om)}  \int_{A_q} \big(g(u_\infty) + \om(q) \big)^{p-1} \,\d \mL^n .
\end{split}
\]
Hence, by noting that $A_q =\big\{g(u_\infty)< 1-2\om(q) \big\}$, which a result of our earlier observations, the above estimate implies
\[
\begin{split}
\nu_p(A_q) & \leq \,  \frac{1}{\mL^n(\Om)} \int_{A_q} \Big(\big(1-2\om(q)\big) + \om(q) \Big)^{p-1} \,\d \mL^n
\\
& \leq \,  \frac{\mL^n(A_q)}{\mL^n(\Om)}  \big(1- \om(q) \big)^{p-1}
\\
& \leq \,  \big(1- \om(q) \big)^{p-1} .
\end{split}
\]
This establishes all the claimed estimates.

To conclude the proof, it remains to establish the claims regarding the limiting case. By standard sequential weak* compactness results in the spaces of Radon measures, together with the boundedness of $(\D u_p)_{p>n/\al}$ in $L^{2\al}(\Om;\R^{N\by n})$ (by virtue of Proposition \ref{proposition9}), we obtain the existence of limit measures $\mu_\infty$, $M_\infty$ and $\nu_\infty$ such that $\mu_p \weakstar \mu_\infty$, $M_p \weakstar M_\infty$ and $\nu_p \weakstar \nu_\infty$ in the corresponding spaces over $\overline{\Om}$, along a subsequence as $p\to \infty$. We now also show the additional properties of the measure $\nu_\infty$. By the weak* lower-semicontinuity of measures on open sets, we have
\[
\nu_\infty(A_q) \,\leq\, \liminf_{p\to \infty} \nu_p(A_q) \, =\, 0,
\]
for any $q$ fixed. Hence, by letting $q\to \infty$, the upper continuity of the measure $\nu_\infty$ implies
\[
\nu_\infty \big( \big\{g(u_\infty)< 1 \big\} \big) \, =\, \lim_{q\to \infty} \nu_\infty \big( \big\{g(u_\infty)< 1-2\om(q)\big\} \big) \, =\, 0.
\]
This implies that $\nu_\infty \big(\{g(u_\infty) =1 \} \big) = \nu_\infty(\overline{\Om})$. Further, since the function $h\equiv 1$ belongs to $C(\overline{\Om})$, the subsequential weak* convergence $\nu_p \weakstar \nu_\infty$ implies
\[
\nu_\infty(\overline{\Om}) \, =\, \int_{\overline{\Om}} h \, \d \nu_\infty \, =\, \lim_{p\to \infty} \int_{\overline{\Om}} h \, \d \nu_p  \, =\, \lim_{p\to \infty} \nu_p(\overline{\Om}) \, \geq \, \lim_{p\to \infty} \frac{1}{1+\om(p)}\, =\,1.
\] 
Finally, since $u_\infty\equiv 0$ on $\p\Om$ and by \eqref{2.2} we get $g(0)=0$, it follows that $g(u_\infty)\equiv 0$ on $\p\Om$. Thus, we infer that $\p\Om \sub \big\{g(u_\infty)<1 \big\}$, which yields that  $\nu_\infty(\p\Om)=0$. \qed

\ms

Now we may establish the satisfaction of the PDE system \eqref{1.5} for the quadruple $(u_\infty,\La_\infty,M_\infty,\nu_\infty)$ under only assumptions \eqref{1.2} and  \eqref{1.3}.

\begin{lemma} \label{lemma11} Let $M_\infty \in \mM(\overline{\Om};\R^{N\by n})$ and $\nu_\infty \in \mM(\overline{\Om})$ be the Radon measures obtained in Lemma \ref{lemma10}. Then, the pair $(u_\infty,\La_\infty)$ satisfies the divergence PDE system \eqref{1.5}, weakly in $(C^1_0(\overline{\Om};\R^N))^*$ (namely \eqref{1.12} holds true for any $\phi \in C^1_0(\overline{\Om};\R^N)$).
\end{lemma}

\BPL \ref{lemma11}. Fix $\phi \in C^1_0(\overline{\Om};\R^N)$ and $p>n/\al+2$.  By \eqref{1.11} and \eqref{1.16}, we may rewrite the divergence PDE system \eqref{1.10} as
\[
\left\{
\begin{array}{rl}
\div \big( \p f(\D u_p) \mu_p \big)  + \, \La_p \,\p g(u_p) \nu_p\,=\,0, & \text{ in }\Om,
\\
u_p \,=\,0, & \text{ on }\p\Om.
\end{array}
\right.
\]
By the measure identity 
\[
M_p \, =\, \p f(\D u_p) \mu_p,
\]
The weak formulation of the PDE means that for any $\phi \in C^1_0(\overline{\Om};\R^N)$, we have
\[
\int_\Om \D \phi : \d M_p   = \, \La_p \int_\Om \p g(u_p) \cdot \phi \,\d \nu_p.
\]
By Proposition \ref{proposition9}, we have that $\La_p\larrow \La_\infty$ and also $u_p \larrow u_\infty$ uniformly on $\overline{\Om}$ as $p\to \infty$ along a subsequence. By \eqref{1.3}, we also have that $\p g(u_p) \larrow \p g(u_\infty)$ uniformly on $\overline{\Om}$ as $p\to \infty$, along the same subsequence. The conclusion follows directly by the application of Lemma \ref{lemma10} and the strong-weak* continuity of the duality pairing $C(\overline{\Om}) \by \mM(\overline{\Om}) \larrow \R$.     \qed

\ms

\section{Regularisations up to the boundary} \label{section5}

In this section we introduce the appropriate mollifications that will be utilised in the next section to show the satisfaction of the equation \eqref{1.14}. This regularisation scheme utilises results on the geometry of (strongly) Lipschitz domain from Hofmann-Mitrea-Taylor \cite{HMT}, and is closely related to the regularisation schemes used in Ern-Guermond \cite{EG}.

To begin with, let $n\in L^\infty(\p\Om,\mH^{n-1};\R^n)$ be the outer unit normal vector field on $\p\Om$. Then, (see Hofmann-Mitrea-Taylor \cite[Sec.\ 2, 4]{HMT} for the proofs of the claims in the paragraph) there exists a vector field $\xi \in C^\infty_c(\R^n ;\R^n)$ that is globally transversal to $n$ on $\p\Om$, namely exists $\de_0>0$ such that
\[
\xi \cdot n \, \geq\, \de_0 , \quad \mH^{n-1}\text{-a.e.\ on }\p\Om.
\]
Further, $\xi$ can be chosen to have length $|\xi|\equiv 1$ in an open collar $\{\dist(\cdot,\p\Om)<r_0\}$ around $\p\Om$ for some $r_0>0$ and to vanish on $\{\dist(\cdot,\p\Om)>2r_0\}$. If $\p\Om$ is a compact $C^\infty$ manifold, then one can choose $\xi:=n$ and the transversality condition is satisfied for $\de_0=1$. Further, there exists $\ell, \e_0>0$ such that, for all $\e \in (0,\e_0)$ we have
\[
\dist\big(x+\e \ell\xi(x),\, \p\Om \big) \geq \, 2\e, \quad\text{ for all }x\in\p\Om.
\]

Using the above observations from \cite{HMT}, we may now define our mollifiers. Fix $\e \in (0,\e_0)$ and some $\varrho \in C^\infty_c(\mB_1(0))$ satisfying $\varrho\geq 0$ and $\|\varrho\|_{L^1(\R^n)}=1$. For any $v\in L^\infty(\Om;\R^N)$, extended by zero on $\R^n\set \Om$, we define $\K^\e v : \R^n \larrow \R^N$ by setting
\[
(\K^\e v)(x) \,:=\, \int_{\R^n}  v\big(x+\e \ell\xi(x)-\e y \big)\, \varrho(y)\, \d y.
\]
What this regularisation does is to ``compress" $v$ to a map which is compactly supported inside $\Om$ before mollifying. However, for technical convenience it is \emph{not} exactly equal to the standard mollifier of the ``compressed" function $v\big(\cdot +\e \ell\xi(\cdot) \big)$, which instead equals the convolution
\[
\frac{1}{\e^n}\varrho\Big(\frac{\cdot}{\e}\Big)* v\big(\cdot +\e \ell\xi(\cdot) \big)
\]
(this would require to put ``$\xi(x-\e y)$" instead of ``$\xi(x)$" in the formula defining $\K^\e$). The advantage of this slight variance, as we will see right next, is a simpler formula for the derivatives.

The next result lists the main properties of this mollification scheme.

\begin{proposition} \label{proposition12} The family of regularisation operators $(\K^\e)_{0<\e<\e_0}$ satisfy the next properties:

\ms

\noi \emph{(1)} For any $v\in L^\infty(\Om;\R^N)$, we have $\K^\e v \in C^\infty_c(\Om;\R^N)$ and also $\K^\e v \larrow v$ as $\e\to 0$ a.e.\ on $\Om$ and in $L^q(\Om;\R^N)$, for any $q\in[1,\infty)$.
\ms

\noi \emph{(2)} For any $v\in L^\infty(\Om; \R^N)$ and any convex function $\Phi :\R^N \larrow \R$ satisfying  $0\in \argmin \{\Phi : \R^N\}$ (namely such that $\Phi \geq \Phi(0)$ on $\R^N$), we have
\[
\Phi\big(\K^\e v(x)\big) \,\leq\, \underset{\Om\cap \mB_\e(x+\e \ell\xi(x))}{\ess\sup}\, \Phi(v) ,
\]
for any $x\in \Om$. In particular, $\Phi(\K^\e v) \leq  \|\Phi(v)\|_{L^\infty(\Om)}$ on $\Om$.
\ms

\noi \emph{(3)} For any $v\in W^{1,\infty}_0(\Om; \R^N)$, we have
\[
\left\{ \ \ 
\begin{split}
\D(\K^\e v) \, & =\, \K^\e (\D v) \big[I + \e\ell\D\xi \big]^\top,
\\
\big| \D(\K^\e v) -\K^\e (\D v) \big| \, &\leq\, \e \ell\|\D\xi \|_{L^\infty(\R^n)} \|\D v \|_{L^\infty(\Om)},
\end{split}
\right.
\] 
on $\Om$. Also, as $\e\to 0$ we have $\K^\e v \larrow v$ in $W^{1,q}_0(\Om;\R^N)$ for all $q\in[1,\infty)$ and in $C^\ga(\overline{\Om};\R^N)$ for all $\ga \in (0,1)$. Further, $\K^\e v \weakstar v$ in $W^{1,\infty}_0(\Om;\R^N)$.

\ms

\noi \emph{(4)}  For any $v\in W^{1,\infty}_0(\Om;\R^N)$ and any convex function $\Phi :\R^{N\by n} \larrow \R$ satisfying $0\in \argmin \{\Phi : \R^{N\by n}\}$, there exists $C>0$ such that 
\[
\Phi\big(\D(\K^\e v)(x)\big) \,\leq\, \underset{\Om\cap \mB_\e(x+\e \ell\xi(x))}{\ess\sup}\, \Phi(\D v) \,+\,\e C^*,
\]
for any $x\in \Om$, where the constant $C^*$ depends only on $\Om$, $\ell$, $\D\xi$, $\D\Phi$ and $\|\D v \|_{L^\infty(\Om)}$. In particular, 
\[
\Phi\big(\D(\K^\e v)\big) \, \leq  \, \|\Phi(\D v)\|_{L^\infty(\Om)} \,+\, \e C^*,  
\]
everywhere on $\Om$.
\end{proposition}

\BPP \ref{proposition12}. (1) A change of variables yields the identity
\[
(\K^\e v)(x) \,=\, \int_{\mB_\e(x+\e \ell\xi(x))} \frac{1}{\e^n}\varrho\bigg(\frac{x+\e \ell\xi(x)- z}{\e}\bigg) \, v(z)\,\d z 
\]
which, combined with the fact that $\dist\big(x+\e \ell\xi(x),\, \p\Om \big) \geq 2\e$ when $x\in \p\Om$, imply that $\K^\e v \equiv 0$ on an open neighbourhood of $\R^n \set \Om$ because $v \equiv 0$ on $\R^n \set \Om$ and also $\smash{\bar \mB_\e(x+\e \ell\xi(x)) \subseteq \R^n \set \overline{\Om}}$  when $x\in \p\Om$ (and also for $x$ in open neighbourhood of $\p\Om$). Further, since $\supp (\varrho) \sub \bar \mB_1(0)$, the integral above is in fact equal to the same integral taken over $\R^n$, hence we easily deduce by recursive differentiation that $\K^\e v \in C^\infty (\R^n;\R^N)$.
\ms

\noi (2) Since $\varrho \mL^n$ is a probability measure on $\R^n$, by Jensen's inequality, we estimate
\[
\begin{split}
\Phi\big(\K^\e v(x)\big) \, & =\, \Phi \left(\,\int_{\R^n}  v\big(x+\e \ell\xi(x)-\e y \big)\, \varrho(y)\, \d y \right)
\\
& \leq\, \int_{\R^n} \Phi \Big(v\big(x+\e \ell\xi(x)-\e y \big)\Big) \varrho(y)\, \d y
\\
&\leq\, \underset{y\in \mB_1(0)}{\ess\sup}\, \Phi \Big(v\big(x+\e \ell\xi(x)-\e y \big)\Big) 
\\
&=\, \underset{\mB_\e(x+\e \ell\xi(x))}{\ess\sup}\, \Phi(v) ,
\end{split}
\]
for any $x\in\Om$. By using that $v \equiv 0$ on $\R^n \set \Om$ and our assumption on $\Phi$, we further have
\[
\begin{split}
\Phi\big(\K^\e v(x)\big) \, &  \leq\, \max\bigg\{ \underset{\Om\cap \mB_\e(x+\e \ell\xi(x))}{\ess\sup}\, \Phi(v),\ \underset{\mB_\e(x+\e \ell\xi(x)) \set\Om}{\ess\sup}\, \Phi(v) \bigg\}
\\
&\leq\, \max\bigg\{ \underset{\Om\cap \mB_\e(x+\e \ell\xi(x))}{\ess\sup}\, \Phi(v),\  \Phi(0) \bigg\}
\\
&\leq\, \underset{\Om\cap \mB_\e(x+\e \ell\xi(x))}{\ess\sup}\, \Phi(v) ,
\end{split}
\]
for any $x\in\Om$.

\ms

\noi (3) We readily compute
\[
\D\big(\K^\e v \big)(x)\, =\, \int_{\R^n}  \D v\big(x+\e \ell\xi(x)-\e y \big) \big[I + \e\ell\D\xi(x) \big]^\top \varrho(y)\, \d y 
\]
which yields the claimed identity. The desired inequality is a simple consequence of the  above identity together with the estimate of Part (2). The asserted modes of convergence follow by standard arguments on mollifiers (see e.g.\ \cite{E}).

\ms

\noi (4) Since $\Phi :\R^{N\by n} \larrow \R$ is convex, it is in $W^{1,\infty}_{\rm loc}(\R^{N\by n})$. Fix $R>0$ such that 
\[
R \,>\, \|\D\xi \|_{L^\infty(\R^n)} \|\D v \|_{L^\infty(\Om)}.
\]
Then, by Parts (2)-(3) we estimate
\[
\begin{split}
\Phi\big(\D(\K^\e v) (x)\big) \, & =\, \Phi\Big(\K^\e (\D v)(x)  \,+\, \e\ell \K^\e (\D v)\D\xi(x)^\top \Big) 
\\
& \leq \, \Phi\big(\K^\e (\D v)(x)\big) \, + \, \|\D \Phi\|_{L^\infty(\mB_R(0))} \big\|\e\ell  \K^\e (\D v)\D\xi^\top \big\|_{L^\infty(\Om)}
\\
& \leq \, \Phi\big(\K^\e (\D v)(x)\big) \, + \, \e\ell \|\D \Phi\|_{L^\infty(\mB_R(0))} \|\D v \|_{L^\infty(\Om)} \|\D\xi \|_{L^\infty(\R^n)} 
\\
&\leq\, \underset{\Om\cap \mB_\e(x+\e \ell\xi(x))}{\ess\sup}\, \Phi(\D v)\, +\, C^*\e,
\end{split}
\]
for any $x\in\Om$, where we have set
\[
C^*\,:=\, \ell \|\D \Phi\|_{L^\infty(\mB_R(0))} \|\D\xi \|_{L^\infty(\R^n)} \|\D v \|_{L^\infty(\Om)}.
\]
The proof of the proposition is now complete.      \qed

\ms

\section{The divergence PDE system in $L^\infty$, part II} \label{section6}

In this section we establish the satisfaction of \eqref{1.14} for the minimising quadruple $(u_\infty,\La_\infty,\mu_\infty,\nu_\infty)$, under the hypotheses \eqref{1.2}, \eqref{1.3} and also \eqref{1.4}. We begin with some notation.

\begin{remark} \label{remark13}
Under \eqref{1.4}, since $f$ is assumed quadratic, \eqref{1.2} is in fact satisfied for $C_1=C_2$ (see the observations following Example \ref{example3}). Hence, by introducing the positive constant 
\[
\ka\,:=\, \frac{C_8}{C_1} \ \bigg( \!\! =\, \frac{C_7}{C_2} \bigg),
\]
the conclusion of Lemma \ref{lemma8} strengthens to the equality
\[
L_p \, =\, \ka^{1/p}\La_p.
\]
\end{remark}

Next, we derive some differential identities and energy inequalities, which will be utilised to obtain the necessary estimates.

\begin{lemma}
\label{lemma14} For any $p \in (n/\al+2,\infty)$, consider the quadruple $(u_p,\La_p,\mu_p,\nu_p)$ as in Section \ref{section4}. Then, for any $v\in W^{1,\infty}_0(\Om;\R^N)$, we have the differential identity
\[
\int_{\overline{\Om}} f(\D v - \D u_p)\,\d\mu_p \, =  \int_{\overline{\Om}} f(\D v) \,\d\mu_p \,- \int_{\overline{\Om}} f(\D u_p)\,\d\mu_p \,+ \, \La_p \int_{\overline{\Om}} \p g(u_p )\cdot (u_p-v)\,\d\nu_p .
\]
\end{lemma}
Note that in view of \eqref{1.11} and \eqref{1.16}, we have that $\nu_p(\p\Om)=\mu_p(\p\Om)=0$ for $p \in (n/\al+2,\infty)$. Thus, all the integrals above are non-trivial only over $\Om$.

\BPL \ref{lemma14}. Since $v-u_p \in W^{1,\al p}_0(\Om;\R^N)$, by using \eqref{1.4} we have
\[
\begin{split}
\int_{\overline{\Om}} f(\D v - \D u_p)\,\d\mu_p \, & = \int_{\overline{\Om}} \A : (\D v - \D u_p) \ot (\D v - \D u_p )\,\d\mu_p,
\end{split}
\]
which, by virtue of Lemma \ref{lemma4}(2), can be expanded as
\[
\begin{split}
\int_{\overline{\Om}} f(\D v - \D u_p)\,\d\mu_p \, &  = \int_{\overline{\Om}} \A : \D v  \ot \D v \,\d\mu_p \, +\, \int_{\overline{\Om}} \A : \D u_p  \ot \D u_p \,\d\mu_p
\\
& \ \ \ + \int_{\overline{\Om}} \A : \Big( \D u_p  \ot  (-\D v)  \,+\, (-\D v)  \ot  \D u_p\Big) \,\d\mu_p
\\
&=\,  \int_{\overline{\Om}} f(\D v) \,\d\mu_p \,-\,\int_{\overline{\Om}} f(\D u_p)\,\d\mu_p 
\\
&\ \ \ + \int_{\overline{\Om}} \A : \Big( \D u_p  \ot (\D u_p -\D v)  \,+\, (\D u_p -\D v) \ot \D u_p \Big) \,\d\mu_p
\\
=\,  \int_{\overline{\Om}} f   &  (\D v) \,\d\mu_p \,-\,\int_{\overline{\Om}} f(\D u_p)\,\d\mu_p \, + \int_{\overline{\Om}} \p f(\D u_p) : \D (u_p - v) \,\d\mu_p.
\end{split}
\]
By arguing as in the proof of Lemma \ref{lemma11} and testing against $\phi:=u_p-v$ in the weak formulation, we readily deduce the claimed identity. \qed

\ms

\begin{lemma}
\label{lemma15} In the setting of Lemma \ref{lemma14}, for any $p \in (n/\al+2,\infty)$ we have
\[
\left\{ \ \ 
\begin{split}
 \int_{\overline{\Om}} f(\D u_p) \,\d\mu_p \, &=\, \ka \La_p,
 \\
  \int_{\overline{\Om}} |\D u_p|^2 \,\d\mu_p \, &\leq \, \frac{\ka}{C_4}\big(\La_p \, +\, C_3 \ka^{-1/p}\big).
\end{split}
  \right.
\]
\end{lemma}

\BPL \ref{lemma15}. From \eqref{1.16}, Lemma \ref{lemma8} and Remark \ref{remark13} (noting also that $\mL^n(\p\Om)=0$), we have
\[
\begin{split}
 \int_{\overline{\Om}} f(\D u_p) \,\d\mu_p \, =\,  \av_{\Om} f(\D u_p) \frac{f(\D u_p)^{p-1}}{\La_p^{p-1}}\,\d\mL^n \, = \, \frac{1}{\La_p^{p-1}}L_p^{p} = \, \ka \La_p.
\end{split}
\]
Further, by \eqref{1.2} and \eqref{1.4} (noting also that $\al=2$ in \eqref{1.2} under \eqref{1.4}), we have
\[
\ka \La_p \, = \int_{\overline{\Om}} f(\D u_p) \,\d\mu_p\, \geq\, C_4 \int_{\overline{\Om}} |\D u_p|^2 \,\d\mu_p \, -\, C_3 \mu_p(\overline{\Om}).
\]
The claimed inequality is now a consequence of the above together with the bound $\mu_p(\overline{\Om}) \leq \ka^{1-1/p}$, which follows from Lemma \ref{lemma10} and Remark \ref{remark13}.  \qed

\ms

The following result is an immediate consequence of assumption \eqref{1.4}, Lemmas \ref{lemma4}, \ref{lemma14}, \ref{lemma15} and Proposition \ref{proposition9}.

\begin{corollary} \label{corollary16} In the setting of Lemma \ref{lemma15}, for any $p \in (n/\al+2,\infty)$ and any map $v\in W^{1,\infty}(\Om;\R^N)$, we have the estimate
\[
c_0 \int_{\overline{\Om}} \big|\D v - \D u_p \big|^2\,\d\mu_p \, \leq \,  \int_{\overline{\Om}} f(\D v) \,\d\mu_p \,- \ka \La_p \,+ \, C_0 \|u_p-v \|_{L^\infty(\Om)} ,
\]
where 
\[
c_0\, :=\, \min\si(\mS \A), \ \ \ \ \ C_0 := \sup_{n/\al+2<p<\infty} \! \big\{\La_p \|\p g(u_p) \|_{L^\infty(\Om)}\big\}.
\]
\end{corollary}

Now we expound on the methodology utilised in the remainder of this section, in order to complete the proof of Theorem \ref{theorem2}.

\begin{remark}[The method] \label{remark17} The estimate of Corollary \ref{corollary16} is the main energy estimate we will need to pass to the limit as $p\to \infty$ along a sequence in the PDE system \eqref{1.10} to obtain \eqref{1.14}. The main difficulty in trying that is that both $\D u_p$ and $\mu_p$ converge in a weak (weak*) sense only, and in fact in different spaces. Therefore, a priori it is not at all clear that
\[
\D u_p \mu_p \weakstar \D u_\infty \mu_\infty
\]
and the product $\D u_\infty \mu_\infty$ may not be a well defined measure, as in general $\D u_\infty$  is in $L^\infty(\Om;\R^{N\by n})$ and only Lebesgue measurable, hence it may not be defined on lower-dimensional subsets of $\overline{\Om}$ on which the Borel measure $\mu_\infty$ may concentrate (e.g.\ on hypersurfaces in $\Om$ or on the boundary $\p\Om$).

To circumvent these problems, we argue as follows. The idea is to show first that 
\[
\D u_p \,\mu_p \weakstar V_\infty \,\mu_\infty, \quad \text{ in }\mM\big(\overline{\Om};\R^{N\by n}\big),
\]
as $p\to \infty$, for some Borel measurable $V_\infty : \overline{\Om} \larrow \R^{N\by n}$ in $L^2(\overline{\Om},\mu_\infty;\R^{N\by n})$. (In particular, the measure $V_\infty \mu_\infty$ is then well defined.) Then, using weak* lower-semicontinuity we let $p\to \infty$ in Corollary \ref{corollary16} to obtain for any map $v \in C^1_0(\overline{\Om};\R^N)$ (not just in $W^{1,\infty}_0({\Om};\R^N)$) that
\[
c_0 \int_{\overline{\Om}} \big| \D v - V_\infty \big|^2\,\d\mu_\infty \, \leq \,  \int_{\overline{\Om}} f(\D v) \,\d\mu_\infty \,- \ka \La_\infty \,+ \, C_0 \|u_\infty-v \|_{L^\infty(\Om)} .
\]
Next, one would like to set ``$v:=u_\infty$" to obtain ``$\D u_\infty = V_\infty$". However, this is not directly possible, at the very least because directly $\D u_\infty \mu_\infty$ is not well defined, and also the limiting process of the previous step prevents us for setting $v=u_\infty$ due to the lack of regularity. To this end, we utilise the regularisation operators $(\K^\e)_{0<\e<\e_0}$ introduced in Section \ref{section5} to set $v:= \K^\e u_\infty$ and deduce as $\e\to 0$ that a \emph{special Borel measurable representative $\D u_\infty^\star$ does exist}, which is a version of the gradient of $u_\infty$ in the equivalence classes of a.e.-equality for both $\D u_\infty \in L^\infty(\Om;\R^{N\by n})$ and for $V_\infty \in L^2(\overline{\Om},\mu_\infty;\R^{N\by n})$. Hence, the PDE system \eqref{1.14} is satisfied for this representative of $\D u_\infty$ which in particular makes $\p f(\D u_\infty^\star)\mu_\infty$ a well-defined measure (recall that $\p f$ is a just linear mapping under \eqref{1.4}).
\end{remark}

Now we deploy the method set out in Remark \ref{remark17}.

\begin{lemma} \label{lemma18} There exists a mapping $V_\infty \in L^2(\overline{\Om},\mu_\infty;\R^{N\by n})$ such that, along a subsequence $(p_j)_1^\infty$ we have 
\[
\D u_p \,\mu_p \, \weakstar \, V_\infty \,\mu_\infty, \quad \text{ in }\mM\big(\overline{\Om};\R^{N\by n}\big),
\]
as $p\to \infty$. Further, for any non-negative continuous function $\Phi \in C\big(\overline{\Om} \by \R^{N\by n}\big)$ such that $X \mapsto \Phi(x,X)$ is convex and of quadratic growth at infinity, we have
\[
\int_{\overline{\Om}} \Phi(\cdot, V_\infty)\, \d\mu_\infty \, \leq\, \liminf_{p\to\infty} \int_{\overline{\Om}} \Phi(\cdot, \D u_p)\, \d\mu_p.
\]
\end{lemma}

\BPL \ref{lemma18}. By Lemma \ref{lemma15} and Proposition \ref{proposition9}, we have that
\[
\sup_{p\in(n/\al+2,\infty)} \int_{\overline{\Om}} |\D u_p|^2\, \d\mu_p \, < \, \infty.
\]
In view of this estimate, the conclusion follows by the theory of measure-function pairs of Hutchinson in \cite[Sec.\ 4, Def.\ 4.1.1, 4.1.2, 4.2.1 and Th.\ 4.4.2]{H}.     \qed

\ms

In virtue of the above considerations, we obtain the following rather immediate consequence.

\begin{lemma} \label{lemma19} In the setting of Lemma \ref{lemma18}, for any $p \in (n/\al+2,\infty)$ and any fixed $v\in C^1_0(\overline{\Om};\R^N)$, we have the estimate
\[
c_0 \int_{\overline{\Om}} \big| \D v - V_\infty \big|^2\,\d\mu_\infty \, \leq \,  \int_{\overline{\Om}} f(\D v) \,\d\mu_\infty \,- \ka \La_\infty \,+ \, C_0 \|u_\infty-v \|_{L^\infty(\Om)} .
\]
\end{lemma}

\BPL \ref{lemma19}. Since by assumption $\D v \in C(\overline{\Om};\R^{N\by n})$, we may apply Lemma \ref{lemma18} to $\Phi(x,X):=|\D v(x)-X|^2$, which satisfies the required convexity and continuity requirements. Then, we use the estimate of Corollary \ref{corollary16}, together with the facts that $(u_p,\La_p) \larrow (u_\infty,\La_\infty)$ in $C(\overline{\Om};\R^N)\by \R$ as $p\to\infty$ along a sequence, and that $f(\D v) \in C(\overline{\Om})$, which is the predual space of $\mM(\overline{\Om})=( C(\overline{\Om}))^*$.   \qed

\ms

The following result is the most remarkable consequence of the energy estimate of Lemma \ref{lemma19}.

\begin{lemma} \label{lemma20}  \emph{(i)} For \emph{any} sequence $(v_j)_1^\infty \sub C^1_0(\overline{\Om};\R^N)$ which satisfies the conditions 
\[
\left\{ \ \ \ \
\begin{split}
 \lim_{j\to\infty}\|v_j - u_\infty \|_{L^\infty(\Om)} \,  &=\, 0, 
\\
 \limsup_{j\to\infty} \| f(\D v_j) \|_{L^\infty(\Om)} \, &\leq\, \La_\infty ,
\end{split}
\right.
\]
we have that
\[
\lim_{j\to\infty} \int_{\overline{\Om}} \big| V_\infty -\D v_j \big|^2\,\d\mu_\infty \, =\, 0.
\]
\emph{(ii)} If additionally to the above, the sequence $(v_j)_1^\infty$ satisfies
\[
 \lim_{j\to\infty}\|\D v_j - \D u_\infty\|_{L^1(\Om)} \, =\, 0, 
\]
then we obtain the additional conclusion
\[
\lim_{j\to\infty} \int_{{\Om}} \big| \D u_\infty -\D v_j \big|^q\,\d \mL^n \, =\, 0,
\]
for any $q\in[1,\infty)$. Hence, for any such sequence $(v_j)_1^\infty$, we have
\[
\left\{ \ \ 
\begin{split}
&\text{$\D v_j \larrow V_\infty$, \ \ \ in $L^2(\overline{\Om},\mu_\infty;\R^{N\by n})$},
\\
&\text{$\D v_j \larrow \D u_\infty$, \ in $L^q({\Om},\mL^n;\R^{N\by n})$},
\end{split}
\right.
\]
as $j\to\infty$, along perhaps a subsequence.
\end{lemma}

We note that, since $L^q(\Om,\R^{N\by n})$ is a uniformly convex space when $q\in (1,\infty)$, the assumption of $\D v_j \larrow \D u_\infty$ in $L^1(\Om;\R^{N\by n})$ can actually be replaced by the weaker condition $\|\D v_j\|_{L^q(\Om)} \larrow \|\D u_\infty\|_{L^q(\Om)}$ for some $q\in (1,\infty)$.

\BPL \ref{lemma20}. (i) It suffices to apply the assumed modes of convergence to the estimate of Lemma \ref{lemma19}, and recall that $\mu_\infty(\overline{\Om})\leq \ka^{1-1/p}$, as a result of Lemma \ref{lemma10} and Remark \ref{remark13}. 
\ms

\noi (ii) By Lemma \ref{lemma4} and \eqref{1.6} (shown in Proposition \ref{proposition9}), the hypothesis 
\[
\| f(\D v_j)\|_{L^\infty(\Om)} \, \leq \, \|f(\D u_\infty)\|_{L^\infty(\Om)} \, +\, o(1)_{j\to\infty}
\]
implies that $(\D v_j)_1^\infty$ is bounded in $L^q(\Om,\R^{N\by n})$, for all $q\in [1,\infty]$. Hence, by passing perhaps to a subsequence, we have that $\D v_j \weak \D u_\infty$
in $L^q(\Om,\R^{N\by n})$ for all $q\in[1,\infty)$ and also $\D v_j \weakstar \D u_\infty$ in $L^\infty(\Om,\R^{N\by n})$. Since $\D v_j \larrow \D u_\infty$ in $L^1(\Om;\R^{N\by n})$, by the Vitaly convergence theorem and the $L^\infty$ gradient bound, it follows that in fact the convergence is strong in $L^q(\Om,\R^{N\by n})$ for all $q\in[1,\infty)$. The result ensues.    \qed

\ms

In the next result we use Lemma \ref{lemma20} identify the limit $V_\infty$ as a version of $\D u_\infty$ after perhaps modification on a Lebesgue nullset (recall also that $\mL^n(\p\Om)=0$).

\begin{corollary} \label{corollary21} There exists a Borel measurable mapping $\D u_\infty^\star : \overline{\Om} \larrow \R^{N\by n}$, which is a version of both $\D u_\infty \in L^\infty(\Om;\R^{N\by n})$ and of $V_\infty \in L^2(\overline{\Om},\mu_\infty;\R^{N\by n})$, namely
\[
\D u_\infty^\star \, =
\left\{
\begin{array}{l}
\D u_\infty, \ \ \mL^n\text{-a.e.\ on }\Om, 
\\
V_\infty, \ \ \ \ \mu_\infty\text{-a.e.\ on }\overline{\Om}.
\end{array}
\right.
\]
Further, $\D u_\infty^\star$ can be represented as
\[
\D u_\infty^\star (x) \, = \left\{
\begin{array}{ll}
\underset{j \to \infty}{\lim}\, \D \big(\K^{\e_j}u_\infty \big)(x), & \text{if the limit exists},
\\
0, & \text{otherwise}.
\end{array}
\right.
\]
along an infinitesimal sequence $(\e_j)_1^\infty \sub (0,1)$, where $(\K^\e)_{0<\e<\e_0}$ are the regularisation operators of Section \ref{section5}.
\end{corollary}

\begin{remark}\label{remark22} It follows that, if one defines the Borel set
\[
G \,:=\, \Big\{ x\in \overline{\Om}\ : \ \not\exists\ \underset{j \to \infty}{\lim}\, \D \big(\K^{\e_j}u_\infty \big)(x) \Big\},
\]
then $G$ is both a $\mL^n$-nullset and a $\mu_\infty$-nullset: $\mL^n(G)=\mu_\infty(G)=0$.
\end{remark}

\BPCOR \ref{corollary21}. Let $(\K^\e)_{0<\e<\e_0}$ be the regularisation operators of Section \ref{section5}. By Proposition \ref{proposition12}, we have $\K^\e u_\infty \in C^\infty_c(\Om;\R^N)$, therefore the choice $v:=\K^\e u_\infty$ in the estimate of Lemma \ref{lemma19} is admissible. Again by Proposition \ref{proposition12} and \eqref{1.6} (proved in Proposition \ref{proposition9}), we have that
\[
\left\{ \ \ \ \
\begin{split}
 \lim_{\e\to 0}\|\K^\e u_\infty - u_\infty \|_{L^\infty(\Om)} \,  &=\, 0, 
\\
 \limsup_{\e\to 0} \big\| f\big(\D (\K^\e u_\infty)\big) \big\|_{L^\infty(\Om)} \, &\leq\, \La_\infty ,
 \\
\lim_{\e\to 0}\|\K^\e u_\infty - u_\infty \|_{W^{1,1}_0(\Om)} \,  &=\, 0.
\end{split}
\right.
\]
We may now apply Lemma \ref{lemma20} to obtain that, for any infinitesimal sequence $(\e_j)_1^\infty$ there is a subsequence, symbolised again by $(\e_j)_1^\infty$, such that for any $q \in [1,\infty)$,
\[
\left\{ \ \ 
\begin{split}
&\text{$\D \big(\K^{\e_j} u_\infty \big) \larrow V_\infty$, \ \ \ in $L^2(\overline{\Om},\mu_\infty;\R^{N\by n})$},
\\
&\text{$\D \big(\K^{\e_j} u_\infty \big) \larrow \D u_\infty$, \ in $L^q({\Om},\mL^n;\R^{N\by n})$},
\end{split}
\right.
\]
as $j\to\infty$. By passing perhaps to a further subsequence, we infer that
\[
\left\{ \ \ 
\begin{split}
&\text{$\D \big(\K^{\e_j} u_\infty \big) \larrow V_\infty$, \ \  \, $\mu_\infty$-a.e.\ on $\overline{\Om}$},
\\
&\text{$\D \big(\K^{\e_j} u_\infty \big) \larrow \D u_\infty$, \ $\mL^n$-a.e.\ on ${\Om}$},
\end{split}
\right.
\]
as $j\to\infty$. The conclusion follows by defining the map $\D u_\infty^\star$ as in the statement. The result therefore ensues.     \qed

\ms

We may now establish the satisfaction of the necessary conditions.

\begin{lemma} \label{lemma23} The quadruple $(u_\infty,\La_\infty,\mu_\infty,\nu_\infty)$ satisfies the system of PDEs \eqref{1.14}, weakly in the dual space $(C^1_0(\overline{\Om};\R^N))^*$.
\end{lemma}

\BPL \ref{lemma23}. By the proof of Lemma \ref{lemma11}, for any fixed $\phi \in C^1_0(\overline{\Om};\R^N)$ we have
\[
\int_\Om \p f(\D u_p):\D \phi \, \d \mu_p   = \, \La_p \int_\Om \p g(u_p) \cdot \phi \,\d \nu_p.
\]
By Proposition \ref{proposition9}, we have $\La_p\larrow \La_\infty$ and also $u_p \larrow u_\infty$ in $C(\overline{\Om};\R^N)$ as $p\to \infty$ along a sequence. Further, $\p g(u_p) \larrow \p g(u_\infty)$ in $C(\overline{\Om};\R^N)$ as $p\to \infty$. By Lemma \ref{lemma18} and Corollary \ref{corollary21}, we have that
\[
\D u_p \,\mu_p \, \weakstar \, \D u_\infty^\star \,\mu_\infty, \quad \text{ in }\mM\big(\overline{\Om};\R^{N\by n}\big).
\]
Further, by virtue of Lemma \ref{lemma4}, $\p f$ is a linear mapping on $\R^{N\by n}$, therefore 
\[
\p f(\D u_p) \,\mu_p \, \weakstar \, \p f(\D u_\infty^\star) \,\mu_\infty, \quad \text{ in }\mM\big(\overline{\Om};\R^{N\by n}\big).
\]
The conclusion follows directly by the application of Lemma \ref{lemma10} and the strong-weak* continuity of the duality pairing $C(\overline{\Om}) \by \mM(\overline{\Om}) \larrow \R$.    \qed
\ms

We conclude with establishing that the set whereon the measure $\mu_\infty$ concentrates is the set whereon $f(\D u^\star_\infty)$ is maximised over $\overline{\Om}$.

\begin{lemma} \label{lemma24} The next equalities hold true:
\[
\left\{ \ \ 
\begin{split}
& \hspace{40pt}  \La_\infty \, =\, \frac{1}{\ka} \int_{\overline{\Om}} \, f(\D u_\infty^\star)\, \d\mu_\infty, 
\\
& \mu_\infty\text{-\,}\underset{\overline{\Om}}{\ess\sup}\, f(\D u_\infty^\star) \, =\,  \underset{\overline{\Om}}{\sup}\, f(\D u_\infty^\star) \, =\, \La_\infty .
\end{split}
\right.
\]
Further, we have
\[
\mu_\infty(\overline{\Om})=\ka, \ \ \ \mu_\infty \left(\big\{f(\D u^\star_\infty) <\La_\infty \big\} \right) =0.
\]
Finally, the boundary $\p\Om$ is a nullset with respect to the Radon measure $\D u_\infty^\star \,\mu_\infty$: 
\[
\big\|\D u_\infty^\star \,\mu_\infty \big\|(\p\Om) \, =\,0.
\]
\end{lemma}

\BPL \ref{lemma24}. By Lemma \ref{lemma19} for $v:=\K^\e u_\infty$, we have
\[
0\, \leq \,  \int_{\overline{\Om}}\, f\big(\D (\K^\e u_\infty)\big) \,\d\mu_\infty \,- \,\ka \La_\infty \,+ \, C_0 \|u_\infty- \K^\e u_\infty \|_{L^\infty(\Om)} .
\]
By Corollary \ref{corollary21}, we have that $\D (\K^\e u_\infty) \larrow \D u_\infty^\star$ in $L^2\big(\overline{\Om},\mu_\infty;\R^{N\by n})$ as $\e_j\to 0$, because $V_\infty=\D u_\infty^\star$ $\mu_\infty$-a.e.\ on $\overline{\Om}$. Further, by Lemma \ref{lemma4}, by Proposition \ref{proposition12} and by the dominated convergence theorem, the above estimate yields as $\e_j\to 0$ that
\[
\La_\infty \, \leq \, \frac{1}{\ka} \int_{\overline{\Om}} \, f( \D u_\infty^\star ) \,\d\mu_\infty  .
\]
On the other hand, recall that by Lemma \ref{lemma10} and Remark \ref{remark13} we have $\mu_\infty(\overline{\Om}) \leq \ka$. Hence, by Proposition \ref{proposition12}, Lemma \ref{lemma20} and H\"older's inequality, we obtain 
\[
\begin{split}
 \La_\infty \, +\, o(1)_{j\to\infty} \, &\geq \,  \sup_{\overline{\Om}}\, f\big(\D (\K^{\e_j} u_\infty)\big)
\\
&\geq \,  \mu_\infty\text{-\,}\underset{\overline{\Om}}{\ess\sup}\, f\big(\D (\K^{\e_j} u_\infty)\big)
\\
&\geq \, \bigg( \frac{1}{\ka}\, \mu_\infty(\overline{\Om})\bigg) \, \mu_\infty\text{-\,}\underset{\overline{\Om}}{\ess\sup}\, f\big(\D (\K^{\e_j} u_\infty)\big)
\\
&\geq\, \frac{1}{\ka} \int_{\overline{\Om}} \,f\big(\D (\K^{\e_j} u_\infty)\big)\, \d\mu_\infty.
\end{split}
\]
By letting $j\to \infty$, this yields
\[
\La_\infty \, \geq \,  \frac{1}{\ka} \int_{\overline{\Om}}\, f( \D u_\infty^\star ) \,\d\mu_\infty  .
\]
The above estimates establish the claimed integral identity. Further, let $G \sub \overline{\Om}$ be the Borel set of Remark \ref{remark22}. Then, by the definition of $\D u_\infty^\star$ as being equal to zero on $G$ and the fact that $f\geq f(0)=0$ as a result of \eqref{2.1}, we have
\[
\begin{split}
\mu_\infty\text{-\,}\underset{\overline{\Om}}{\ess\sup}\, f(\D u_\infty^\star) \, &=\, \underset{\overline{\Om}  \set G}{\sup}\, f(\D u_\infty^\star)
\\
&=\, \max\bigg\{ \underset{\overline{\Om}  \set G}{\sup}\, f(\D u_\infty^\star) ,\ f(0)\bigg\}
\\
&=\, \max\bigg\{ \underset{\overline{\Om}  \set G}{\sup}\, f(\D u_\infty^\star) ,\ \underset{G}{\sup}\, f(\D u_\infty^\star)\bigg\}
\\
&=\, \underset{\overline{\Om}}{\sup}\, f(\D u_\infty^\star).
\end{split}
\]
Arguing similarly for the Lebesgue measure, we obtain 
\[
\begin{split}
\La_\infty \, =\, \| f(\D u_\infty^\star) \|_{L^\infty(\Om)}\, =\, \underset{\overline{\Om}  \set G}{\sup}\, f(\D u_\infty^\star)
=\, \underset{\overline{\Om}}{\sup}\, f(\D u_\infty^\star).
\end{split}
\]
Therefore, the desired equalities have been established. Further, since $\mu_\infty(\overline{\Om}) \leq \ka$ and $f(\D u_\infty^\star) \leq \La_\infty$ on $\overline{\Om}$, we deduce also that in fact $f(\D u_\infty^\star)$ equals its supremum $\La_\infty$ over $\overline{\Om}$ and also the measure of $ \mu_\infty(\overline{\Om})$ is full, namely $ \mu_\infty(\overline{\Om})= \ka$. Finally, by Lemma \ref{lemma19} for $v=:\K^\e u_\infty$ and by using that $V_\infty=\D u_\infty^\star$ $\mu_\infty$-a.e.\ on $\overline{\Om}$ and that $\supp(\K^\e v) \Subset \Om$, we conclude
\[
\begin{split}
c_0 \int_{\p\Om} \big| \D u_\infty^\star \big|^2\,\d\mu_\infty \, & \leq\, c_0 \int_{\overline{\Om}} \big| \D u_\infty^\star -\D \big(\K^{\e_j} u_\infty\big) \big|^2\,\d\mu_\infty 
\\ 
&\leq \,  \int_{\overline{\Om}} f\big(\D \big(\K^{\e_j} u_\infty\big)\big) \,\d\mu_\infty \,- \ka \La_\infty \,+ \, C_0 \|u_\infty- \K^\e u_\infty \|_{L^\infty(\Om)} 
\\
&\leq\, o(1)_{j\to\infty}.
\end{split}
\]
This implies that
\[
\big\|\D u_\infty^\star \,\mu_\infty \big\|(\p\Om) \, \leq \, \sqrt{\mu_\infty(\p\Om)}\bigg(\int_{\p\Om} \big| \D u_\infty^\star \big|^2\,\d\mu_\infty\bigg)^{\!\!1/2} =\, 0.
\]
The result ensues.
\qed
\ms

\section{Computations for an explicit example} \label{section7}

In this section we provide some explicit computations in the case of the $\infty$-eigenvalue problem on the ball. Let $n,N\geq 1$ and choose $f:=\frac{1}{2}|\cdot|^2$
 and $g:=\frac{1}{2}|\cdot|^2$ the corresponding Euclidean norms on $\R^{N\by n}$ and $\R^N$ respectively. We also fix $R>0$ and choose $\Om := \mB_R$, the ball of radius $R$ centred at the origin of $\R^n$. Then, for any direction $e\in \R^N$ with $|e|=1$, the directed cone function
\[
 C_\infty(x) := \Big(1 -\frac{|x|}{R}\Big)e  
\]
is a solution to 
\[
\| \D C_\infty  \|_{L^\infty(\mB_R)}  = \, \inf \Big\{ \| \D v  \|_{L^\infty(\mB_R)}  : \ v  \in W^{1,\infty}_0(\mB_R;\R^N), \ \| v \|_{L^\infty(\mB_R)} =1\Big\}.
\]
Indeed, we have $C_\infty \in W^{1,\infty}_0(\mB_R;\R^N)$ and $\| C_\infty \|_{L^\infty(\mB_R)} =1$, because $0 \leq |C_\infty| \leq 1$ and also $|C_\infty (0)|=1$. Further, for any $x\in \mB_R\set \{0\}$ we have $|\D C_\infty(x)| = \frac{1}{R}$, which yields 
\[
\big\| \D C_\infty \big\|_{L^\infty(\mB_R)}=\frac{1}{R}. 
\]
This value in fact is the infimum over all maps in the admissible class. To see this, fix any $v  \in W^{1,\infty}_0(\mB_R;\R^N)$ with $\| v \|_{L^\infty(\mB_R)} =1$. By continuity, this means that $0\leq |v| \leq 1$ and that exists $\bar x \in \mB_R$ such that $|u(\bar x)|=1$. Let $(\K^\e)_{0<\e<\e_0}$ be the family of regularising operators introduced earlier in Section \ref{section5}. Then, $\K_\e v \in C^1_c(\mB_R;\R^N)$ and by standard arguments on viscosity solutions (see e.g.\ \cite[Section 3]{K0}) there exists $\bar x_\e \in \mB_R$ close to $\bar x$ such that $|\K^\e v|$ attains its maximum at $\bar x_\e$ and $|\K^\e v(\bar x_\e)| \larrow 1$ as $\e \to 0$. Then, for any fixed $z \in \p \mB_\R$ we have $\K^\e v(z)=0$ and hence
\[
\begin{split}
1+o(1)_{\e\to 0} \,&= \,\big|\K^\e v(\bar x_\e) - \K^\e v(z) \big| 
\\
& = \, \left|\int_0^1 \D (\K^\e v) \big(\la \bar x_\e +(1-\la) z \big) \cdot (\bar x_\e -z)\,\d \la\right|
\\
& \leq\, \| \D v\|_{L^\infty(\mB_R)} |\bar x_\e -z|.
\end{split}
\]
By choosing $z:=R \frac{\bar x_\e}{|\bar x_\e|}$ if $\bar x_\e \neq 0$ and any $z \in \p \mB_R$ if $\bar x_\e=0$, we obtain by letting $\e \to 0$ that 
\[
\big\| \D v \big\|_{L^\infty(\mB_R)} \,\geq \, \frac{1}{R}. 
\]
The above arguments show that the directed cone $C_\infty$ is indeed a vectorial $\infty$-eigenfunction. 

We now consider the necessary PDEs that $C_\infty$ solves. By invoking Theorems \ref{theorem1}-\ref{theorem2} and by noting that $\La_\infty = \big\| \D C_\infty \big\|_{L^\infty(\mB_R)}= 1/R$, we see that there exist measures $\mu_\infty,\nu_\infty \in \smash{\mM(\bar \mB_R)}$ such that
\[
-\div \big(\D C_\infty \mu_\infty \big)  \, =\, \frac{1}{R}C_\infty \nu_\infty \ \text{ \ in }\Om.
\]
Since $\nu_\infty(\bar \mB_R)=1$ and in this case $\{C_\infty = \sup_{\mB_R}C_\infty\}=\{0\}$, we deduce that 
\[
\nu_\infty \, =\, \de_0.
\]
Therefore, since $C_\infty\de_0 = C_\infty(0)\de_0$, the PDE system reduces to
\[
-\div \big(\D C_\infty \mu_\infty \big)  \, =\, \frac{1}{R}e\de_0 \ \text{ \ in }\Om.
\]
Further, since $\D C_\infty = -\frac{1}{R} e \ot \sgn$, where $\sgn$ is the sign function in $\R^N$, we may compute an explicit measure $\mu_\infty$, which in fact is absolutely continuous  on $\bar \mB_R$. By using the fundamental solutions of the Laplacian $\De$
\[
\Phi(x) = -\frac{1}{2\pi}\ln|x| \  \text{ for }n=2,\ \ \  \Phi(x) = \frac{1}{n(n-2)\al(n) |x|^{n-2}} \ \text{ for }n\geq 3
\]
(where $\al(n)$ symbolises the volume of the unit ball in $\R^n$), which in both cases give 
\[
\D \Phi(x) \, =\, -\frac{1}{n \al(n) |x|^{n-1}}\sgn(x)
\]
and that the system reduces to the single PDE 
\[
-\div( -\sgn \mu_\infty ) \,=\,\de_0 \ \text{ \ in }\Om,
\]
we obtain that the PDE is satisfied for the absolutely continuous measure
\[
\mu_\infty \, =\, \frac{1}{n\al(n)|\cdot|^{n-1}}\mL^n \LL_{\mB_R}.
\]
\begin{remark} \label{remark25} It is worth noting that, as shown in \cite{BMa}, due to the full rotational symmetry and the homogeneity of the vectorial $p$-eigenvalue problem, all vectorial $p$-eigenfunctions are essentially scalar. Even though this is not automatically true for the $\infty$-eigenvalue problem, it does carry over to $p=\infty$ at least for those $\infty$-eigenfunctions which are constructed as $L^p$-limits. Nevertheless, this reduction to essentially scalar minimisers is not deducible for the problem \eqref{1.1} for general $f,g$.
\end{remark}

\subsection*{Acknowledgments} The author would like to thank Erik Lindgren for various scientific discussions and his expert insights on the $\infty$-eigenvalue problem, which took place during an academic visit of the author to the University of Uppsala in late 2019, as well as in various electronic communications thereafter.


\bibliographystyle{amsplain}

\end{document}